\documentclass[11pt]{article}
\usepackage{amsmath}
\usepackage{amssymb}
\usepackage{float}
\usepackage{graphicx}

\usepackage[comma,sort&compress]{natbib}
\usepackage{color} 

\usepackage{lineno}
\usepackage{url}
\usepackage{pdfpages}
\usepackage{setspace} 
\usepackage{lineno}

\usepackage{hyperref}

\usepackage[table]{xcolor}

\hypersetup{
	bookmarks=true,         % show bookmarks bar?
	unicode=false,          % non-Latin characters in Acrobat‚????s bookmarks
	pdftoolbar=true,        % show Acrobat‚????s toolbar?
	pdfmenubar=true,        % show Acrobat‚????s menu?
	pdffitwindow=false,     % window fit to page when opened
	pdfstartview={FitH},    % fits the width of the page to the window
	pdftitle={My title},    % title
	pdfauthor={Author},     % author
	pdfsubject={Subject},   % subject of the document
	pdfcreator={Creator},   % creator of the document
	pdfproducer={Producer}, % producer of the document
	pdfkeywords={keywords}, % list of keywords
	pdfnewwindow=true,      % links in new window
	colorlinks=true,       % false: boxed links; true: colored links
	linkcolor=blue,          % color of internal links
	citecolor=blue,        % color of links to bibliography
	filecolor=magenta,      % color of file links
	urlcolor=cyan           % color of external links
}

\usepackage[labelfont=bf,labelsep=period,justification=raggedright]{caption}

\makeatletter
\renewcommand{\@biblabel}[1]{\quad#1.}
\makeatother

\date{}

\title{\bf{%Incentives for cooperative agreement compliance: enhancing participation is as important as discouraging misbehaviour / Incentive for commitment compliance in PD / centralised commitments / 
Institutional Incentives for the Evolution of Committed Cooperation:  Ensuring  Participation is as Important as Enhancing Compliance} }

\author{The Anh Han$^{1,\star}$}

\begin{document}
	\maketitle
	{\footnotesize
		\noindent
		$^{1}$ School of Computing, Engineering and Digital Technologies,  Teesside University, Middlesbrough, UK TS1 3BA\\
		$^\star$ Corresponding: The Anh Han (t.han@tees.ac.uk) 
	}

\newpage
\section*{Abstract}

 Both conventional wisdom and empirical evidence suggest that arranging a prior commitment or agreement
before an interaction takes place enhances the chance of reaching mutual cooperation. Yet it is not clear  what mechanisms might underlie  the participation in and compliance with such  a commitment, especially when participation is costly and non-compliance can be  profitable. Here, we develop a theory of participation  and compliance with respect to an explicit commitment formation process and to institutional incentives where individuals, at first, decide  whether or not to join a cooperative agreement to play a one-shot social dilemma game. Using a mathematical model, we determine whether and when participating in a costly commitment, and complying with it, is an evolutionary stable strategy,  resulting in high levels of cooperation. We show that, given a sufficient budget for providing incentives, rewarding of commitment compliant behaviours better promotes cooperation than punishment of non-compliant ones. Moreover, by sparing part of this  budget for rewarding those willing to participate in a  commitment,  the overall level  of cooperation can be  significantly enhanced for both reward and punishment. Finally,  the presence of mistakes in deciding to  participate  favours evolutionary stability of commitment compliance  and cooperation.

 \noindent \textbf{Keywords:} Commitment, reward, punishment,  evolution of cooperation, social dilemma, evolutionary dynamics. 

 \newpage
 
%\doublespacing

\section{Introduction}

Commitments, such as contracts and agreements, are fundamental components of many social and economic interactions, ranging from personal, to institutional, to political or religious ones, in order to ensure a mutually beneficial  outcome for the parties involved \citep{IronsChapter2001,nesse2001evolution,HanBook2013,akdeniz2021evolution,frank88,cherry2013enforcing,sasaki2012take}. %They can be  in the form of formal  contracts, informal social norms or  non-binding promises \citep{shelton2003commitment,nesse2001evolution}. 
Arranging a prior commitment from all parties involved before an interaction  improves the chance that people can reach mutual cooperation when individual interests are in conflict %, as shown in public good behavioural experiments
\citep{chen1994effects,kerr1997still,balliet2010communication}.  
In most modern societies,  institutions are created to enforce formal contracts and enhance cooperation, through suitable incentive structures such as punishment for wrongdoing  \citep{zumbansen2007law,ostrom2009understanding,nesse2001evolution}. 
People joining a religion share certain norms and expectations and might expect certain reward and punishment for a given behaviour \citep{johnson2006hand,IronsChapter2001}. 
Commitments are important  in the context of computerised multi-agent systems, where they are formalised as a tool for regulating agents' interactions and collective behaviour \citep{singh2013norms}.
%Historically, communities have constructed institutions that provide incentives and enforcement such as reward and punishment, in order to ensure compliance to  commitments \citep{ostrom1990governing,nesse2001evolution}. 

Evolutionary game theory (EGT) \citep{key:Sigmund_selfishnes} provides an appropriate tool to study the evolution of  cooperative behaviour in  social dilemmas, as they are  governed by institutional incentives \citep{sasaki2012take,sigmundinstitutions,chen2015first,wang2019exploring,DuongHanPROCsA2021,CIMPEANU2021107545,gois2019reward,sun2021combination} and prior  commitments  \citep{han2013good,HanJaamas2016,sasaki2015commitment,akdeniz2021evolution}. 
 However, prior works have not examined the interplay between different forms of incentive and commitment behaviours, including participation in and compliance with a prior commitment. 
On the one hand, existing institutional incentive models  do not   capture  the commitment formation process explicitly  %\citep{krapohl2021instability,SantosPNAS2011,sasaki2012take,sigmundinstitutions}.
%Moreover, where   incentives are concerned,  existing research  only target behaviours in the interaction 
\citep{sasaki2012take,sigmundinstitutions,chen2015first,wang2019exploring,gois2019reward}. These works fail  to consider  the need for improving participation in a commitment before the interaction. 
Indeed, since participating in a commitment usually involves a cost, for example initial time and effort spent setting up the commitment \citep{tappin2015financial} and/or  membership fees \citep{heidar2006party}, and requires the involved parties to follow certain restrictive terms and conditions, such participation might need to be encouraged.   
%Indeed, providing incentives for joining a commitment in the first place, as are the case for  such as climate change agreements and  safety agreements,  and  religions,  is quite common (CITES).} 
Examples of incentives for encouraging  participation are many, and  have been shown to be crucial for ensuring positive outcome, as  in the contexts of climate change agreements \citep{barrett2003increasing}, healthcare programs for  reducing smoking  during pregnancy \citep{tappin2015financial} and diabetes \citep{bruni2009economic}. Lottery was used as a form of reward for participation in Covid-19 vaccination in many countries including the US  \citep{sehgal2021impact}.

On the other hand, existing models of commitment formation  assume that a compensation   from a commitment violating player is transferred internally to a commitment compliant co-player,  via institutional enforcement   \citep{Han:2014tl,han2013good,martinez2015apology,hanTom2016synergy}. 
%Since individuals can decide whether not to honour an adopted commitment---with abundant evidence  of commitment breaching in both controlled experiments and real-world scenarios \citep{nesse2001evolution,dannenberg2016non,kerr1997still,nguyen2019}---
These works have thus been limited  in  considering what alternative mechanisms can ensure high levels of commitment compliance and cooperation. For instance, as will be analysed below, one might ask whether positive (reward) or negative (punishment) incentive is more efficient, to that end.  Moreover, similar to the aforementioned institutional incentive models, these works fail to consider the possibility that enhancing participation in a commitment before interaction can maximize the benefit provided by commitments for promoting the evolution of cooperation.

Here, we aim to fill this gap by  investigating, theoretically, how to distribute a per capita incentive budget effectively between enhancing the level of participation in a commitment and ensuring compliance with any adopted commitment, thereby optimising  the overall level of cooperation. Our analysis is carried out in the context of the one-shot Prisoners' Dilemma (PD) game, in which players can either cooperate (C) and defect (D). Before a PD game, players can choose whether or not to join a commitment and cooperate in the game. The commitment is formed if both players  commit, and the players can decide  how to act in the game conditionally on whether such a commitment is formed or not.  In the former case, the committed players share a participation cost ({denoted by $\epsilon$, see also Table \ref{table:parameters}}) as a fee for maintaining the institution that provides incentives. % to enforce compliance to the agreed upon behaviour. 

%We will examine which type of incentive is more efficient for enhancing commitment compliance and cooperation. We study when commitment compliant behaviours are evolutionarily viable, by considering competition among all possible strategies. Furthermore, given the cost of participation, we hypothesise that, if a fraction of the incentive budget is used to encourage players to join the commitment, the overall cooperation outcome can be improved.   %We consider both evolutionary stable strategies (ESS) 
%A strategy is defined by three decisions: i) accepts (A) or not (N) to join the commitment; ii) cooperates (C) or defects (D) if the commitment is formed; iii) cooperates (C) or defects (D) if it is not formed. Thus, there are eight possible strategies in total, which are summarized in  Table \ref{table:eight-strategies}. 

%We also consider different types of noise ...  

The present analysis is based on two well-adopted, complementary approaches in EGT, namely, evolutionary stable strategies (ESS analysis) \citep{maynard-smith:1982to,Otto2007AEvolution} and finite population dynamics \citep{nowak:2005:nature,key:Sigmund_selfishnes}. While the former allows a simple assessment of when a strategy can resist invasion from all other strategies in a population (i.e. being an ESS), it does not capture the detailed stochastic dynamics  among all strategies in co-presence. This drawback is overcome with the latter approach.  We  determine when participating in a costly commitment, and complying with it, is an ESS and promotes the evolution of enhanced cooperation.  
We  examine which incentive, reward of commitment compliant behaviours or punishment of non-compliant ones, is more efficient to this end, and when. 
This analysis  sheds light on whether the cooperation outcome can be enhanced by sparing part of the incentive  budget for encouraging participation in a commitment despite reducing the incentive budget available for promoting  compliant behaviours. 
Finally, we  examine the impact of noise in a participation decision on the  stability of commitment compliance  and  cooperation.

\section{Models and Methods }
\subsection{Model}
We consider a well-mixed population of individuals who, before an interaction, can choose whether or not to voluntarily join a prior commitment dictating that, if it is formed, they then must cooperate in the interaction. Interactions are modeled using  the one-shot Prisoner's Dilemma (PD) game with a payoff matrix of the form 
\begin{equation}
 \bordermatrix{~ & C & D\cr
                  C & R,R & S,T \cr
                  D & T,S & P,P  \cr
                 }.
\end{equation}
If both players choose C (D), they both receive the same reward $R$ (penalty $P$) for mutual cooperation (mutual defection).  Unilateral cooperation provides the sucker's payoff $S$ for the cooperative player and the temptation to defect $T$ for the defecting one.  The payoff matrix corresponds to the preferences associated with the PD when the parameters satisfy the ordering, $T > R > P > S$  \citep{coombs1973reparameterization}.  
To provide  a clearer interpretation of  analytical conditions presented below,  we sometimes reduce the PD to its special case, the donation game \citep{key:Sigmund_selfishnes}, where $T = b,\ R = b-c,\ P = 0,\ S = -c$, with $b$ and $c$ stand for  the benefit and cost of cooperation, respectively.

%The commitment stands when all parties agree. 
%Otherwise, players just interact using the regular PD, in absence of a commitment.  In the former case, the committed players share a cost $\epsilon$. 
\begin{table}
\begin{tabular}{ p{2cm}|p{2.8cm}p{4.5cm}p{4cm}  }
 \hline
Strategies & Accept commitment? & Cooperate in presence of commitment? & Cooperate in absence of commitment? \\
 \hline
 ACC   & Yes    &  Yes  &  Yes \\
 ACD&   Yes  & Yes   &No\\
 ADC &Yes & No&  Yes\\
 ADD    &Yes & No&  No\\
 NCC&   No  & Yes&Yes\\
 NCD& No  & Yes   &No\\
 NDC& No  & No&Yes\\
  NDD& No  & No&No\\
 \hline
\end{tabular}
\caption{The eight strategies with commitment/agreement formation.}
\label{table:eight-strategies}
\end{table}
We  perform a comprehensive analysis where the full set of strategies in considered. Namely, given that players can choose i) whether to accept  (A) or not (N) to join a prior commitment before a PD game, ii) to cooperate (C) or defect (D) in the PD if the commitment is formed, and iii) to cooperate (C) or defect (D) in the PD if the commitment is not formed, in total we can define eight possible strategies. They are  denoted as ACC, ACD, ADC, ADD, NCC, NCD, NDC and NDD, as summarised in Table \ref{table:eight-strategies}. {Also, recall that  a commitment is formed when both players in a PD commit and in that case   the committed players share a participation cost  $\epsilon$.}

We assume that there is a per capita budget $u$ available for providing incentives. A fraction  of the budget, $\alpha u$ ($0 \leq \alpha \leq 1$), is used for rewarding those who are  willing to participate in a commitment (i.e.  AXY players for $X, Y \in \{C,D\}$), increasing the chance a commitment being formed. The remaining budget, i.e. $(1-\alpha) u$, is used for rewarding commitment compliant players (i.e. ACC and ACD players) or punishing non-compliant ones (i.e. ADC and ADD players). When $\alpha = 0$, it means the budget is used only for incentivising commitment compliant  behaviours (i.e., \textit{pure reward} and \textit{pure punishment} scenarios). 
As we consider reward and punishment separately, without loss of generality, we assume that all the incentives described above are equally cost efficient, where the incentive recipient's increased or decreased amount (corresponding to reward and punishment, respectively) equals the institution's cost.

Finally, consider that, with some small probability ${\chi}$, an error might occur when players decide whether to join a commitment. It can be due to fuzzy mind or {trembling} hands, as well as miscommunication or environmental noise.  We show that this type of noise strongly influences  evolutionary dynamics, even in favor of cooperation and commitment compliance. In Supporting Information (\textbf{SI}), we also study other types of noise that occur in an PD interaction.  %The derivation of all the relevant payoff matrices in presence of noise is described in Methods.
\begin{table}%[h]
\centering
\begingroup
\renewcommand{\arraystretch}{1.25}
\begin{tabular}{lc}
\hline
 Parameter  &  Symbol %& Range/Value Analysed
 \\ 
\hline
 Population size & $N$ %& 100 
 \\
  Cost of cooperation & $c$ %& 100 
 \\
 Benefit of cooperation & $b$ %& 100 
 \\
 Intensity of selection & $\beta$  %& \{$0.1$\}  
 \\
% Payoff matrix  & R, S, T, P  & \{$1, -1, 2, 0$\} \\
 Per capita budget & $u$  %& $[0,6]$  
 \\
 Fraction of the budget for rewarding participation  & $\alpha$  %& $[0,1]$  
 \\
 Cost of commitment participation  & $\epsilon$  %&  $[0,6]$   
 \\
Error probability in participating in a commitment  & ${\chi}$  %& $[0,0.2]$ 
\\
\end{tabular}
\endgroup
 \caption{Model parameters }
 \label{table:parameters}
\end{table}
  
The model parameters are summerised in Table \ref{table:parameters}. Next, the derivation of  payoff matrices for all the scenarios described above is {provided}.
%\subsection*{Payoff derivation}
\paragraph{In absence of participation errors.}
First, when incentives are not in use, i.e. the \textit{no policy} scenario {(that is, when $u = 0$)}, the payoff matrix for the eight strategies  (see Table \ref{table:eight-strategies}), reads (for row player) %\footnote{We consider that players might make errors when deciding whether to join an agreement, with a probability $p$. } 
{\footnotesize
\begin{equation} 
 \label{eq:payoff_matrix_no_policy}
\Pi_0=\bordermatrix{~ & ACC & ACD &  ADC &  ADD &  NCC &  NCD &  NDC &  NDD \cr
            ACC & R-\epsilon/2 & R-\epsilon/2 & S-\epsilon/2 & S-\epsilon/2 & R & S & R & S  \cr
            ACD & R-\epsilon/2 & R-\epsilon/2 & S-\epsilon/2 & S-\epsilon/2 & T & P  & T & P \cr
            ADC & T-\epsilon/2 & T-\epsilon/2 & P-\epsilon/2 & P-\epsilon/2 & R & S & R &  S  \cr
           ADD & T-\epsilon/2 & T-\epsilon/2 & P-\epsilon/2 & P-\epsilon/2 & T & P & T & P   \cr
           NCC & R & S & R & S & R & S & R & S  \cr
        NCD  & T & P & T & P & T & P & T & P  \cr
            NDC  & R & S & R & S & R & S & R & S  \cr
            NDD  & T & P & T & P & T & P &T& P  \cr
                  %FREE  & 0 & R & T & T & P & P & P & P & P \cr
                  }.
\end{equation}}\\
An observation is that NDD and NCD are equivalent, so are NCC and NDC, because an agreement is only formed when both players agree to join. Thus, when one of these strategists are involved in an interaction, a commitment is never formed and only the move in absence of a commitment matters. % 
Moreover, ACC and ACD are neutral, so are ADC and ADD. 

%Now, assume that there is a per capita budget $u$ available for providing incentives. It can be  used to reward those who honour an agreement, or punish those dishonour it. Moreover, we examine if cooperation can be improved by by spending a fraction  of the budget, $\alpha u$ ($0 \leq \alpha \leq 1$), for rewarding those who are  willing to participate in the agreement (regardless of whether the agreement is formed), as it might improve the chance a cooperation agreement is formed \footnote{Our additional results show that if reward for participation is only made when the agreement is formed, then this act is always detrimental for cooperation; the larger the fraction used for this purpose, the worse the outcome is}.%, while the remaining $(1-\alpha) u$, is for rewarding good behaviour or punishing bad one during the interaction  (when an agreement is formed). 
 %Note that when $\alpha = 0$, it means the budget is used only for incentivising commitment-based behaviours (i.e., \textit{pure reward and punishment scenario}). 
%For simplicity, we assume that all types of incentives are equally cost efficient, and that the amount spent by the institution equals the amount obtained the the incentive recipient.  

When a per capita budget $u$ is available to \textit{reward commitment compliant behaviours}, with a fraction $\alpha$ of it being used for rewarding participation, the payoff matrix reads
{\footnotesize
\begin{equation} 
 \label{eq:payoff_matrix_reward}
\Pi_R=\bordermatrix{~ & ACC & ACD &  ADC &  ADD &  NCC &  NCD &  NDC &  NDD \cr
            ACC & R-\epsilon/2 + u & R-\epsilon/2 + u  & S-\epsilon/2 + u  & S-\epsilon/2 + u  & R+ \alpha u & S+ \alpha u & R+ \alpha u & S + \alpha u \cr
            ACD & R-\epsilon/2+ u  & R-\epsilon/2 + u & S-\epsilon/2 + u & S-\epsilon/2 + u & T+ \alpha u & P+ \alpha u  & T+ \alpha u & P+ \alpha u \cr
            ADC & T-\epsilon/2 + \alpha u & T-\epsilon/2+ \alpha u & P-\epsilon/2+ \alpha u & P-\epsilon/2+ \alpha u & R+ \alpha u & S + \alpha u& R+ \alpha u &  S+ \alpha u  \cr
           ADD & T-\epsilon/2 + \alpha u& T-\epsilon/2 + \alpha u& P-\epsilon/2+ \alpha u  & P-\epsilon/2 + \alpha u& T+ \alpha u & P+ \alpha u & T+ \alpha u & P + \alpha u  \cr
           NCC & R & S & R & S & R & S & R & S  \cr
        NCD  & T & P & T & P & T & P & T & P  \cr
            NDC  & R & S & R & S & R & S & R & S  \cr
            NDD  & T & P & T & P & T & P &T& P  \cr
                  %FREE  & 0 & R & T & T & P & P & P & P & P \cr
                  }.
\end{equation}}
When a per capita budget $u$ is available to \textit{punish commitment non-compliant behaviours}, with a fraction $\alpha$ of it being used for rewarding participation reads
{\small
\begin{equation} 
 \label{eq:payoff_matrix_punishment}
\Pi_P=\bordermatrix{~ & ACC & ACD &  ADC &  ADD &  NCC &  NCD &  NDC &  NDD \cr
            ACC & \gamma_3 & \gamma_3  & \gamma_4  & \gamma_4 & R+ \alpha u & S+ \alpha u & R+ \alpha u & S + \alpha u \cr
            ACD & \gamma_3  & \gamma_3 & \gamma_4 & \gamma_4 & T+ \alpha u & P+ \alpha u  & T+ \alpha u & P+ \alpha u \cr
            ADC & \gamma_1 & \gamma_1  &\gamma_2 &\gamma_2 & R+ \alpha u & S + \alpha u& R+ \alpha u &  S+ \alpha u  \cr
           ADD & \gamma_1 & \gamma_1 & \gamma_2& \gamma_2& T+ \alpha u & P+ \alpha u & T+ \alpha u & P + \alpha u  \cr
           NCC & R & S & R & S & R & S & R & S  \cr
        NCD  & T & P & T & P & T & P & T & P  \cr
            NDC  & R & S & R & S & R & S & R & S  \cr
            NDD  & T & P & T & P & T & P &T& P  \cr
                  %FREE  & 0 & R & T & T & P & P & P & P & P \cr
                  },
\end{equation}}
where we denote $\gamma_1 = T-\epsilon/2 + (2\alpha-1)u$, $\gamma_2 = P-\epsilon/2 + (2\alpha-1)u$, $\gamma_3 =R-\epsilon/2 + \alpha u$, and  $\gamma_4 =S-\epsilon/2 + \alpha u$, just for the purpose of a neat presentation. 
\paragraph{In presence of participation errors.}
We assume that with a small probability ${\chi}$,  players made an error in the decision whether or not to join an agreement in the commitment formation stage (e.g. due to fuzzy mind or trembling hands).
All the payoff matrices above can be re-written as follows.
Denote $\bar{A} = N$, $\bar{N} = A$, $\bar{C} = D$ and $\bar{D} = C$. 
For $\Pi \in \{\Pi_0,\Pi_R, \Pi_P\}; \ P_j \in \{A,N\};\ X_j, Y_j \in \{C,D\}$,
the payoff when a player $P_1X_1Y_1$ against $P_2X_2Y_2$, $\Pi_{P_1X_1Y_1,P_2X_2Y_2}$,  can be written as
$$
(1-{\chi})^2 \Pi_{P_1X_1Y_1,P_2X_2Y_2} + {\chi}(1-{\chi}) \left(\Pi_{\Bar{P_1}X_1Y_1,P_2X_2Y_2} + \Pi_{P_1X_1Y_1,\bar{P_2}X_2Y_2} \right) + {\chi}^2 \Pi_{\bar{P_1}X_1Y_1,\bar{P_2}X_2Y_2}.
$$ 

\subsection{Evolutionary Stable Strategies (ESS)}
As common assumptions in ESS analysis \citep{Otto2007AEvolution}, we assume that i) mutations are rare and thus, there is at most one mutant strategy $m$ at a time in a population of individuals with resident strategy $r$, and ii) the mutant's effect is negligible on the dynamics. To know if a strategy can be invaded or not by another, we need to compute the difference of absolute fitness between a mutant strategy in a population of resident strategy. If the fitness of the mutant is greater than that of  the resident, the mutant invades the population and becomes resident. If the fitness of the mutant is lower, the mutant disappears and the resident resists invasion. When the two values of fitness are equal, the resident also resists invasion because in an infinitely large population, a mutant strategy cannot invade by drift.
A strategy is ESS if it resists invasion from all other strategies.

\subsection{Evolutionary Dynamics in Finite Population}
In finite population  settings,  individuals'  payoff represents their \emph{fitness} or social \emph{success}, and  evolutionary dynamics is shaped  by social learning \citep{key:Hofbauer1998,key:Sigmund_selfishnes}, whereby the  most successful agents will tend to be imitated more often by the other agents. In the current work, social learning is modeled using  the so-called pairwise comparison rule \citep{traulsen2006}, a  standard approach in EGT,  assuming  that an agent $A$ with fitness $f_A$ adopts the strategy of another agent $B$ with fitness $f_B$ with probability $p$ given by the Fermi function, 
$p_{A, B}=\left(1 + e^{-\beta(f_B-f_A)}\right)^{-1}.$
The parameter  $\beta$ represents  the `imitation strength' or `intensity of selection', i.e., how strongly the agents  base their decision to imitate on fitness difference between themselves and the opponents. For $\beta=0$,  we obtain the limit of neutral drift -- the imitation decision is random. For large $\beta$, imitation becomes increasingly deterministic.
 In line with previous works and human behavioural experiments \citep{Szabo2007,zisisSciRep2015,randUltimatum}, we set $\beta = 0.1$ in the main text.
 
% {Despite using small mutation, it's been shown that the results are robust for larger mutation rates. Also, the analytical conditions do not depend on the small mutation limit assumption. }
 
In the absence of mutations or exploration, the end states of evolution are inevitably monomorphic: once such a state is reached, it cannot be escaped through imitation. We thus further assume that, with a certain mutation probability,  an agent switches randomly to a different strategy without imitating  another agent.  In the limit of small mutation rates, the dynamics will proceed with, at most, two strategies in the population, such that the behavioral dynamics can be conveniently described by a Markov Chain, where each state represents a monomorphic population, whereas the transition probabilities are given by the fixation probability of a single mutant \citep{key:imhof2005,key:novaknature2004}. The resulting Markov Chain has a stationary distribution, which characterizes the average time the population spends in each of these monomorphic end states. 

Let $N$ be the size of the population. Denote $\pi_{X,Y}$ the payoff a strategist X obtains in a pairwise interaction with strategist $Y$ (defined in the payoff matrices). Suppose there are at most two strategies in the population, say, $k$ agents using strategy A ($0 \leq k \leq N$)  and $(N-k)$ agents using strategies B. Thus, the (average) payoff of the agent that uses  A and B can be written as follows, respectively, 
\begin{equation} 
\label{eq:PayoffA}
\begin{split} 
\Pi_A(k) &=\frac{(k-1)\pi_{A,A} + (N-k)\pi_{A,B}}{N-1},\\
\Pi_B(k) &=\frac{k\pi_{B,A} + (N-k-1)\pi_{B,B}}{N-1}.
\end{split}
\end{equation}

Now, the probability to change the number $k$ of agents using strategy A by $\pm$ one in each time step can be written as\citep{traulsen2006} 
\begin{equation} 
T^{\pm}(k) = \frac{N-k}{N} \frac{k}{N} \left[1 + e^{\mp\beta[\Pi_A(k) - \Pi_B(k)]}\right]^{-1}.
\end{equation}
The fixation probability of a single mutant with a strategy A in a population of $(N-1)$ agents using B is given by \citep{traulsen2006,key:novaknature2004}
\begin{equation} 
\label{eq:fixprob} 
\rho_{B,A} = \left(1 + \sum_{i = 1}^{N-1} \prod_{j = 1}^i \frac{T^-(j)}{T^+(j)}\right)^{-1}.
\end{equation} 
%In the limit of neutral selection (i.e. $\beta = 0$), $\rho_{B,A}$ equals the  inverse of population size, $1/N$. 

Considering a set  $\{1,...,q\}$ of different strategies, these fixation probabilities determine a transition matrix $M = \{T_{ij}\}_{i,j = 1}^q$, with $T_{ij, j \neq i} = \rho_{ji}/(q-1)$ and  $T_{ii} = 1 - \sum^{q}_{j = 1, j \neq i} T_{ij}$, of a Markov Chain. The normalized eigenvector associated with the eigenvalue 1 of the transposed of $M$ provides the stationary distribution described above \citep{key:imhof2005}, describing the relative time the population spends adopting each of the strategies.

\paragraph{Risk-dominance.} An important measure to compare the two strategies A and B is which direction the transition is stronger or more probable, an A mutant  fixating in a population of agents using B, $\rho_{B,A}$,  or a B mutant fixating in the population of agents using A, $\rho_{A,B}$. It can be shown that the former is stronger, in the limit of large $N$,  if \citep{key:novaknature2004,key:Sigmund_selfishnes}  
%\begin{equation} 
%(N-2)\pi_{A,A} + N\pi_{A,B} > (N-2)\pi_{B,A} + N\pi_{B,B}
%\end{equation} 
%which,is simplified to 
\begin{equation} 
\label{eq:compare_fixprob_cond_largeN}
\pi_{A,A} + \pi_{A,B} > \pi_{B,A} +  \pi_{B,B}.
\end{equation}
%Furthermore, the fraction $\frac{\rho_{B,A}}{\rho_{A,B}}$ is proportional to the difference of the left hand side of the inequality \footnote{The fraction between the two transition probabilities is an important prediction for the success of A in a pairwise comparison with B. For instance, if we consider a population with two strategies A and B, then the frequency of A is an increasing function of the fraction \citep{nowak:2006bo}.  }. 

\section{Results}

%{\rowcolors{2}{green!80!yellow!50}{green!70!yellow!40}

%}
%There are eight such strategies in total. Among them, there are six strategies that play D either in the presence or in the absence of an agreement. 

%ACD, ADC, ADD, 

%NCD, NDC, NDD 

%We examine systematically how to punish these defectors influence the cooperative outcome. Indeed, we first consider punishing each type separately and compare the effects. Then, we consider punishing increasing larger subsets of these defective types. 

%We can categorize these defective types into 1) two groups: accepting defectors (A**) and refusing defectors (N**); 2) three groups: defecting only in absence of agreement (**D), defecting only in presence of agreement (*D*), and always defecting (*DD). 
%\subsection{Payoff matrix}

  \begin{figure}
    \centering
    \includegraphics[width=0.7\linewidth]{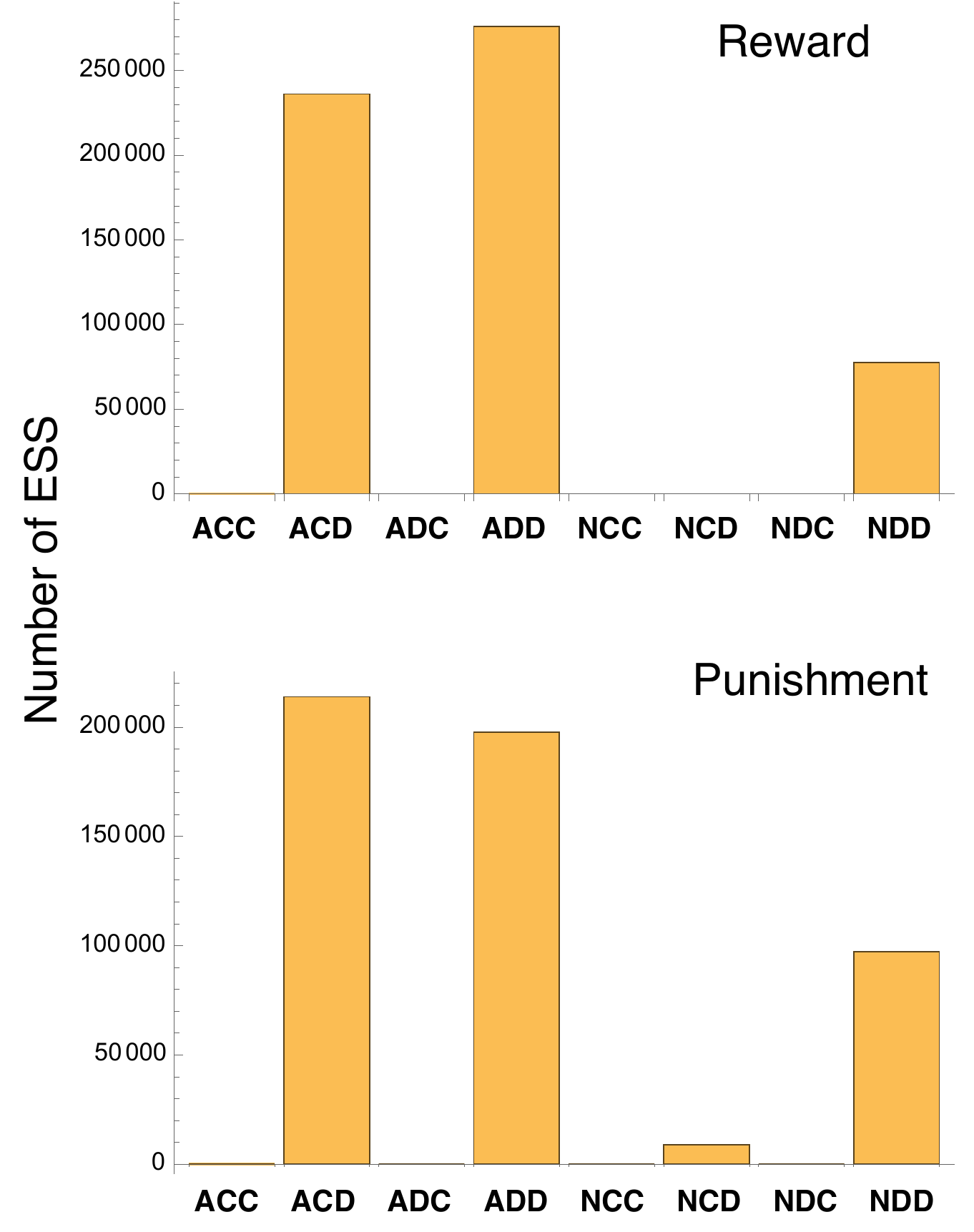}
    \caption{\textbf{Both commitment compliant and non-compliant strategies can be ESS}. We show which strategies can be ESS and their frequencies across the parameters space: $u \in [0,3]$ (increment 0.05), $\epsilon \in [0,3]$ (increment 0.05), $\alpha \in [0,1]$ (increment 0.05) and ${\chi} \in [0,0.2]$ (increment 0.02), i.e. the total number of configurations is thus $61 \times 61 \times 21 \times 11 = 859551$.)   We show the number of times (if any) each strategy is an ESS. We observe that for reward, three strategies ACD, ADD and NDD can be ESS, while  in case of punishment, NCD can also be an ESS.  %When $u$ is small (less then $c = 1$), only defector ESS (ADD and NDD) are possible. When $u > 1$, ACD are ESS frequently. 
    Other parameters:   $R = 1, \ S = -1, \ T = 2, \ P = 0$.  }
    \label{fig:ESS_analysis_overall}
\end{figure}

  \begin{figure}
    \centering
    \includegraphics[width=0.9\linewidth]{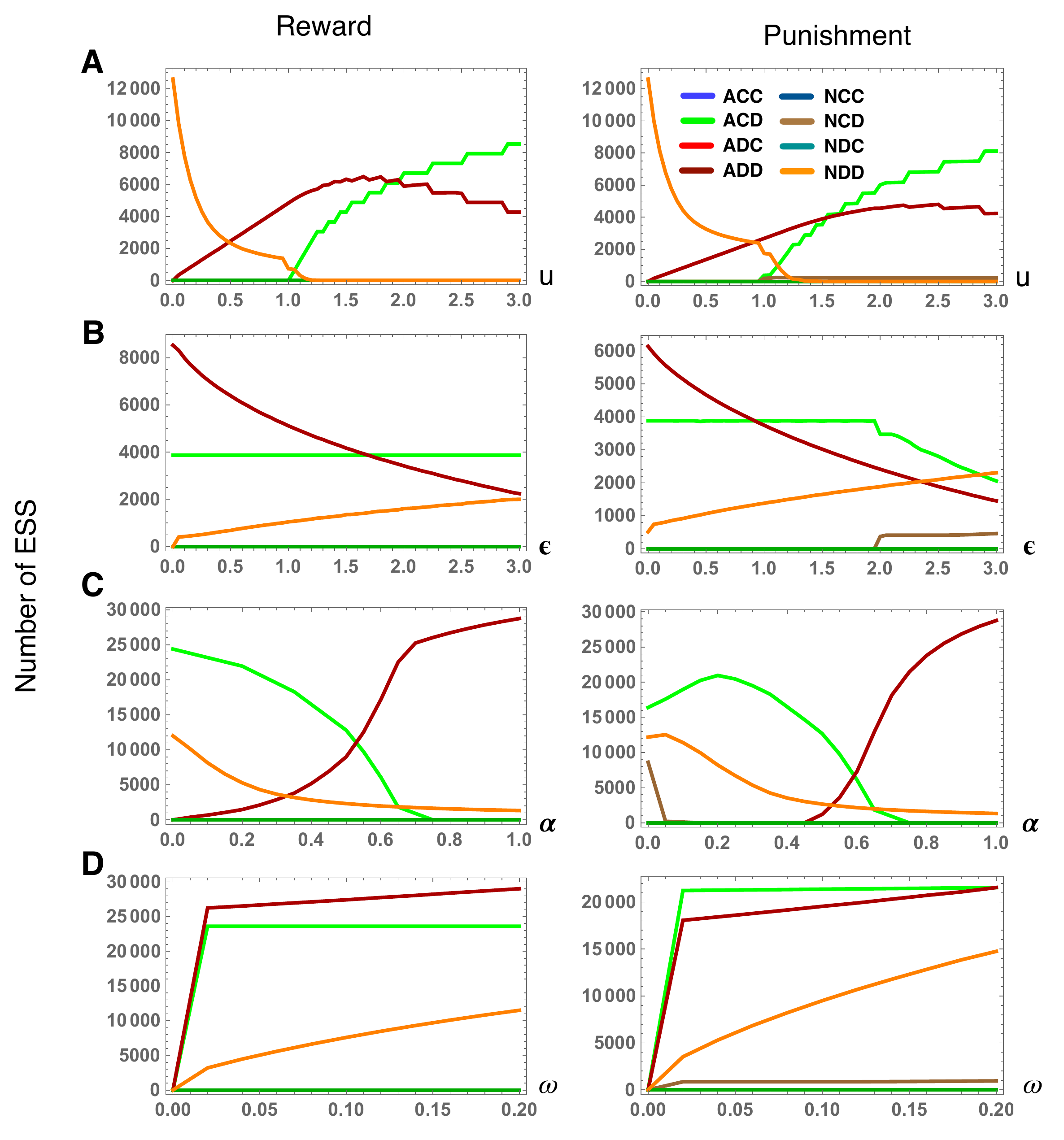}
    \caption{\textbf{For commitment compliance to be evolutionarily stable, it is necessary that noise be non-negligible and  the incentive budget be
 sufficiently high}. {Considering the parameters space  sampled as in Figure 1, i.e. $u \in [0,3]$ (increment 0.05), $\epsilon \in [0,3]$ (increment 0.05), $\alpha \in [0,1]$ (increment 0.05) and ${\chi} \in [0,0.2]$ (increment 0.02), we count  the number  of times each strategy is ESS} for varying $u$ (panel A), $\epsilon$ (panel B), $\alpha$ (panel C) and ${\chi}$ (panel D). We observe that for reward, three strategies ACD, ADD and NDD can be ESS, while  in case of punishment, NCD can also be an ESS. In particular, in both cases, commitment compliant behaviour (ACD) can be ESS only if  $u$ is sufficiently large, $\alpha$ is not too large, and ${\chi}$ is greater than zero.   %When $u$ is small (less then $c = 1$), only defector ESS (ADD and NDD) are possible. When $u > 1$, ACD are ESS frequently. 
    Other parameters:   $R = 1, \ S = -1, \ T = 2, \ P = 0$.  }
    \label{fig:ESS_analysis}
\end{figure}

\begin{figure}
    \centering
    \includegraphics[width=\linewidth]{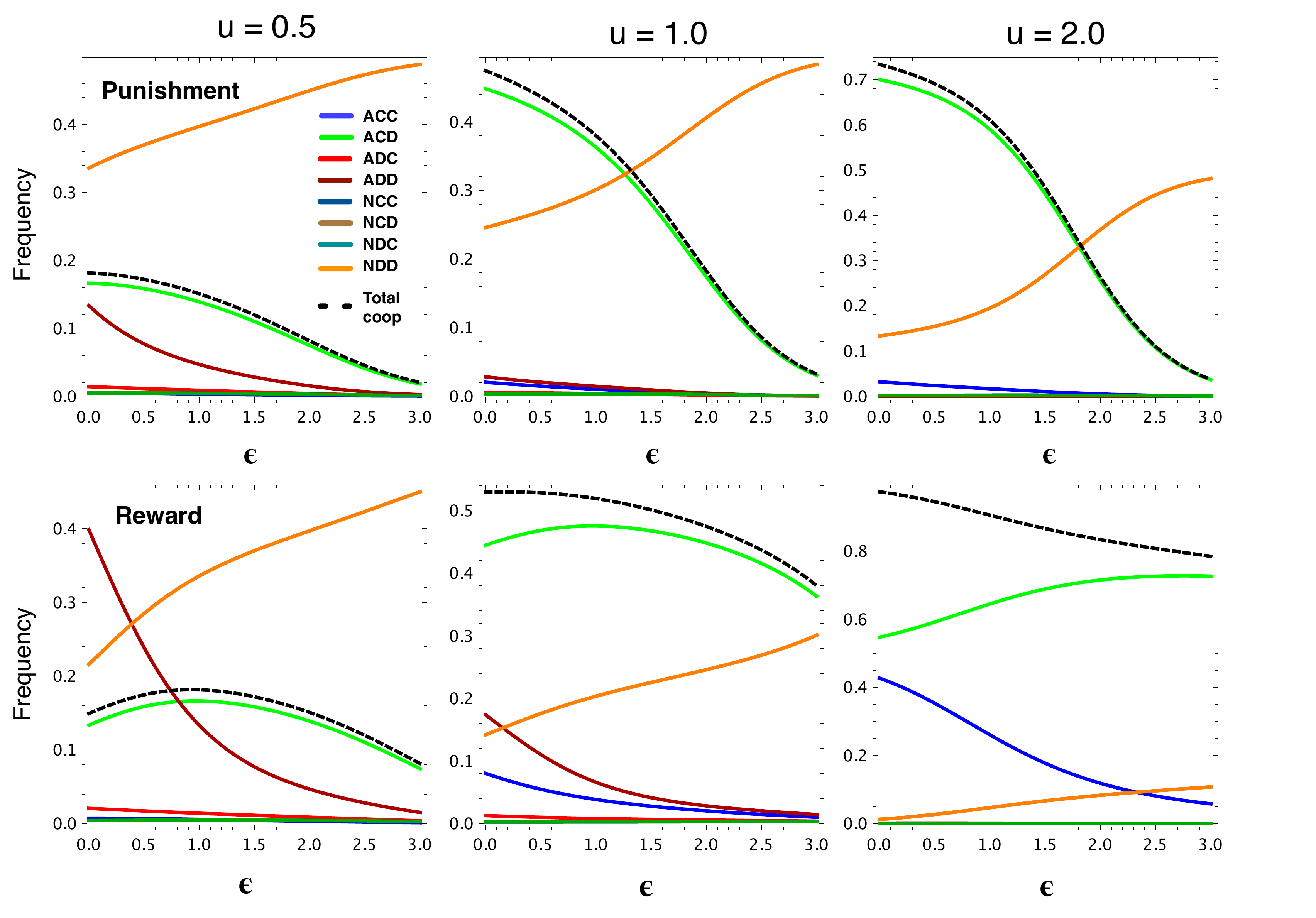}
    \caption{\textbf{Pure incentives promote high levels of commitment compliance and cooperation.} Depicted is the frequency of strategies when pure reward or punishment ($\alpha = 0$) is applied, for different values of the per capita budget ($u$). Both reward and punishment allow ACD to prevail when $\epsilon$ is small. ACC is  more frequent in case of reward {than punishment}. ACD is frequent for a larger range of $\epsilon$ in case of reward  than punishment. These together lead to higher levels of cooperation in case of reward than punishment, for the whole range of $\epsilon$. Reward helps better suppress non-committers when $\epsilon$ is high, explaining its success.   Other parameters: population size $N = 100$,   $R = 1, \ S = -1, \ T = 2, \ P = 0$.   }
    \label{fig:strategies-frequency-no-participation-reward}
\end{figure}

\subsection{ESS Analysis: when complying with a cooperative commitment can be an evolutionarily stable strategy}

%{ADD SOME ANALYTICAL CONDITIONS}
In Figure \ref{fig:ESS_analysis_overall}, we study which strategies can be ESS (see Methods) across the parameter space, namely for $u \in [0,3]$, $\epsilon \in [0,3]$, $\alpha \in [0,1]$ and ${\chi} \in [0,0.2]$.  We observe that, for reward, three strategies ACD, ADD, NDD can be ESS, while NCD can also be ESS in the case of punishment. ACD is the only  ESS that leads to an overall  cooperative outcome (i.e. \textit{cooperative ESS}), where individuals choose to accept a prior commitment and comply with it. All other possible  ESS lead to a defective outcome. Namely,  ADD commits to cooperate but then dishonors the commitment and defects in the interaction, and both NDD and NCD refuse to commit and defect in the interaction.  
%No polymorphism is observed (configuration where two strategies are ESS simultaneously).  

Figure \ref{fig:ESS_analysis} shows the number of times  each strategy is an ESS, for varying $u$, $\epsilon$, $\alpha$ and ${\chi}$ separately. We observe that, in the absence of errors when deciding whether to join a commitment (i.e. ${\chi} = 0$), none of the strategies can be ESS (Figure \ref{fig:ESS_analysis}D), for both types of incentive. For any strategy, there is always another mutant strategy {to which  it is neutral to in absence of errors; namely, for AXY, it is AXY$^\prime$, while for NXY, it is NX$^\prime$Y (where $X \neq X^\prime$,  $Y \neq Y^\prime \in \{C,D\}$). Intuitively,  in the presence of a commitment, the decision of a commitment strategy outside a commitment does not matter, while in the absence of a commitment, the decision of a non-committing strategy inside a commitment does not matter too.}
When  error is non-negligible (${\chi} > 0$), ACD, ADD and NDD can be ESS in case of reward, while NCD can also be ESS in case of punishment. 

Focusing on ACD, it can be an ESS only when a sufficient per capita budget $u$ is available (covering at least the cost of cooperation, $c = 1$, see \textbf{SI} for analytical proof), for both types of incentive (see Figure \ref{fig:ESS_analysis}A). It can be ESS for the whole range of $\epsilon$ being considered{, i.e. up to 3.0} (Figure \ref{fig:ESS_analysis}B). {A noteworthy difference between the two types of incentive is that, when $\epsilon$ is sufficiently large (approximately from 2.0 onward), the number of times ACD being ESS decreases in case of punishment, while it remains stable in case of reward. It implies that reward might outperform  punishment for a larger range of $\epsilon$.}  Moreover, a sufficient fraction  of the budget (i.e. not too large $\alpha$) needs to be spent for rewarding  the compliant  strategies or sanctioning the non-compliant ones (see Figure \ref{fig:ESS_analysis}C). While the frequency of ACD to be ESS decreases with $\alpha$ for reward, it peaks at some intermediate value of $\alpha$ for punishment. This suggests that it might be more important  to incentivise participation for punishment than for reward, in order to increase the chance of commitment compliant behaviour (ACD) to be evolutionarily stable. 

In \textbf{SI} (Figures S1 and  S2), we show that the observations for reward do not change for varying the benefit to cost ratio $b/c$. For punishment, the frequency of ACD to be ESS increases with  $b/c${, reaching a similar frequency as in case of reward when $b/c$ is sufficiently high (around $b/c = 3$). It implies that for commitment compliant behaviours to be evolutionarily stable, it requires the social dilemma  interaction to be more beneficial in case of punishment than reward}. The frequency for NCD decreases with  $b/c$, so that NCD is no longer an ESS when this ratio is sufficiently high. 

Overall, the ESS analysis provides us with some initial insights regarding when commitment compliance (ACD) can be an evolutionarily viable strategy. However, it does not show the detailed dynamics of the whole system and is limited in quantitative characterisation, e.g. of the  overall cooperation level in the population for a given  parameters' configuration. Also, as studied below using a stochastic evolutionary dynamics approach,  ACD and some other strategies can be the most frequent strategy in the population even when they are not  ESS, for example, when noise is absent.  
%To complement this weakness of the ESS analysis, subsequently we study the  stationary distribution of eight strategies adopting a finite population, stochastic approach (CITES). 

\begin{figure}
    \centering
    \includegraphics[width=\linewidth]{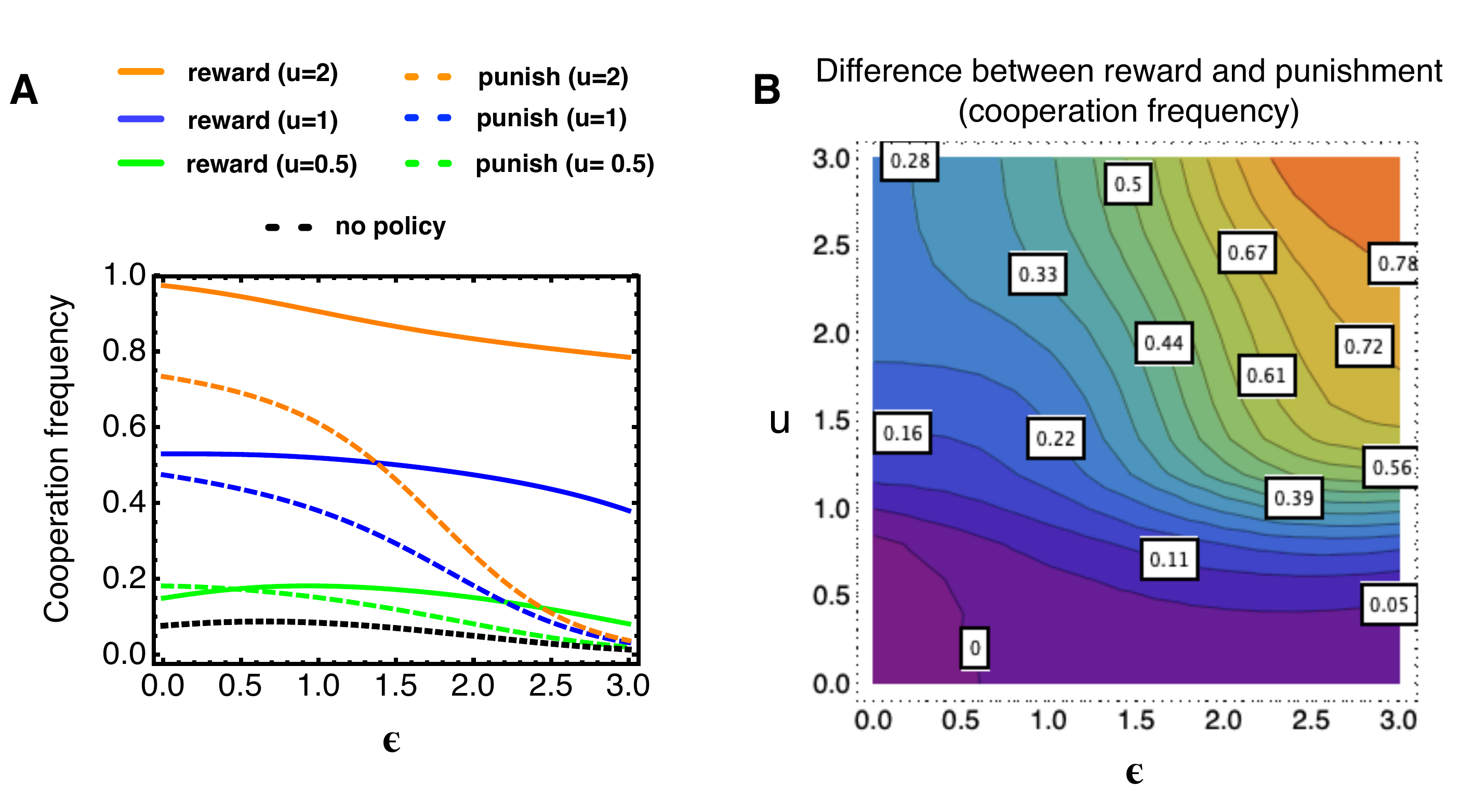}
    \caption{\textbf{Reward leads to a higher level of cooperation than punishment provided a sufficient per capita budget $u$ (the same for both incentives).}  (\textbf{A}) Frequency of cooperation as a function of $\epsilon$  when either reward of commitment compliant strategies (ACC and ACD)  or punishment of non-compliant ones (ADC and ADD), is applied. The black dotted line corresponds to a reference scenario where no policy {(i.e. $u = 0$)} is applied. To ensure a high frequency of cooperation, a sufficient budget for providing   incentives ($u$) is required. Cooperation is reduced when $\epsilon$ increases, for both reward and punishment. Nevertheless, reward leads to a higher level of cooperation than punishment in most cases, except for when both $u$ and $\epsilon$ are sufficiently low (see $u = 0.5$ and $\epsilon \lesssim 0.5$).  (\textbf{B}) Contour plot showing the difference between cooperation obtained through reward and punishment, for varying  $u$ and $\epsilon$. The larger $u$ and $\epsilon$ are, the greater the difference is. Other parameters: population size $N = 100$,  $R = 1, \ S = -1, \ T = 2, \ P = 0$, $\alpha = 0$ (no reward for participation).  }
    \label{fig:cooperation-level-no-reward-participation}
\end{figure}

\subsection{Reward vs punishment for promoting frequent   committed cooperation}
We comparatively study the capability of institutional reward and punishment for promoting the evolution of commitment compliant behaviour and cooperation. 
We first focus on clarifying the effects of pure incentives, considering $\alpha = 0$. We  then study the effect of varying $\alpha$, i.e. when part of the per capita incentive budget is used to reward participation in a commitment before an interaction. We also study the impact of having some small non-negligible error probability  (${\chi} > 0$).

\subsubsection{Pure incentives}
In Figure \ref{fig:strategies-frequency-no-participation-reward}, we compute the long-term frequency (i.e. stationary distribution) of the eight strategies (see Methods) under the pure reward and punishment policies, for varying $\epsilon$ and  different values of $u$. 
Both reward and punishment can enable ACD to prevail when $\epsilon$ is small, given a sufficient budget $u$ ($u = 1$ and 2).  ACC is more frequent in case of reward {than punishment}.  ACD is  frequent and even dominates the population, for a larger range of $\epsilon$ in case of reward {than punishment}. As such, a higher level of overall population cooperation is achieved in case of  reward, for the whole range of $\epsilon$.
When $u$ is small ($u = 0.5$), either ADD or NDD dominates the population, leading to the dominance of defection. NDD also dominates even when $u$ is larger if $\epsilon$ is high, in case of punishment. 
%Reward helps suppress non-committers when $\epsilon$ is high, explaining its success (linked to RISK DOM analysis). 

In \textbf{SI}, {we show the threshold of $\epsilon$ for which ACD is risk-dominant (see Methods) against all other strategies, except for ACC to which it is neutral}. For $\alpha = 0$, for a sufficient budget (namely, $u > (T + P - R - S)/2 = 1$), the thresholds  are
$\epsilon < 2(u +  \min\{T-S, R - P\})$ and 
$\epsilon < 2 \min\{T-S, R - P\}$, respectively, for reward and punishment.  Thus, reward allows for a larger range of $\epsilon$ for which ACD is an evolutionarily viable strategy, having a high long-term frequency. Intuitively, reward enables a better suppression of non-committing strategies such as NDD because the latter do not suffer punishment since they do not commit, while reward provides a payoff advantage for committing strategies {against non-committing ones}.   Notably, the thresholds in numerical results are in accordance with these theoretical observations.  For example, when $u = 2.0$, the threshold is $\epsilon < 2$ for punishment and $\epsilon < 6$ for reward. For $u = 1.0$, since ACD is also neutral to defective committers (ADC and ADD), ACD  is most frequent in the population up to  a slightly slower threshold than the theoretical one (i.e. $\epsilon < 2$ for punishment $\epsilon < 4$ for reward).   
%Note that since ACD is neutral to ACC so the actual threshold might be slightly smaller than the one obtained from the risk dominant analysis. 
%; given that in this case ACD is neutral against defective committers, we can expect a lower threshold than these. For $u = 2.0$, the thresholds is $\epsilon < 2$ for punishment $\epsilon < 6$ for reward. ACD is also risk dominant against defective committers so as expected, the actual  thresholds where ACD dominates the population is very close to the theoretically predicted through risk dominance (although still note that ACD is always neutral to ACC).

In Figure \ref{fig:cooperation-level-no-reward-participation}, we compare the total cooperation frequency in the population obtained through applying pure reward and pure punishment, and also when  no policy  is applied. 
As can be seen, reward leads to a higher level of cooperation than punishment in most cases, except for when both $u$ and $\epsilon$ are rather low.  Also, we observe that the larger $u$ and $\epsilon$ are, the more efficient reward is compared to punishment in terms of cooperation promotion (see Figure \ref{fig:cooperation-level-no-reward-participation}B; also Figure S3 in \textbf{SI}). Figures S9 and S10 in \textbf{SI} show that these observation is robust for varying the intensity of selection $\beta$.
%It is the  total frequency of ACC, ACD, NDC and NCC in the stationary distribution. As observed in Figure 1, we can observe that ACD contributed most to this total cooperation, which is in line with previous works on regimented/enforceable commitments (Han et al, Sci Rep 2013; MORE). 

\begin{figure}
    \centering
    \includegraphics[width=\linewidth]{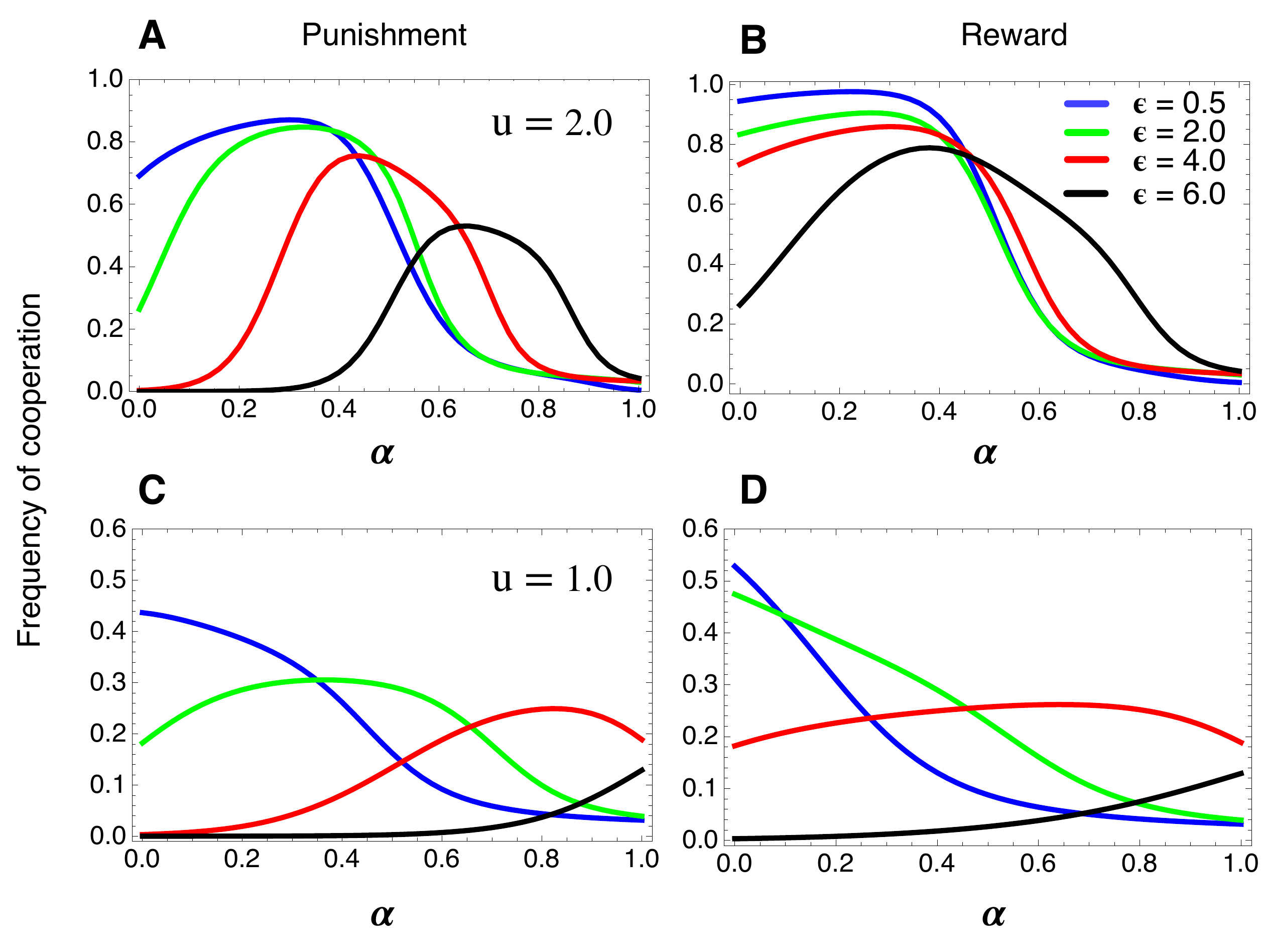}
    \caption{\textbf{Rewarding participation can improve cooperation despite reducing the budget for incentivising commitment compliant behaviour}. Depicted is the frequency of cooperation as a function of the fraction of the budget for rewarding of participation ($\alpha$), for different values of the cost of commitment participation $\epsilon$. Punishment of non-compliant strategies (panels A, C) or reward of complaint ones (panels B, D), are applied as before, using the remaining budget after rewarding participation. We consider scenarios with a large ($u = 2.0$, top row, panels A, B) and small ($u = 1.0$, bottom row, panels C, D) per capita budget. When $\alpha = 0$, it reproduces the results for pure punishment and reward. For both types of incentive, rewarding participation can improve the overall cooperation, especially for the larger $u$. For a larger $\epsilon$, a larger fraction $\alpha$ of the budget should be used for rewarding participation to reach an optimal level of cooperation. %Only when both $u$ and $\epsilon$ are sufficiently small, it is better off not rewarding participation (e.g. when $\epsilon = 0.5$ and $u = 1.0$). 
    Other parameters: $N = 100$,  $R = 1, \ S = -1, \ T = 2, \ P = 0$.   }
    \label{fig:reward-vs-punishment-with-reward-participation-cooperation}
\end{figure}

\begin{figure}
    \centering
    \includegraphics[width=\linewidth]{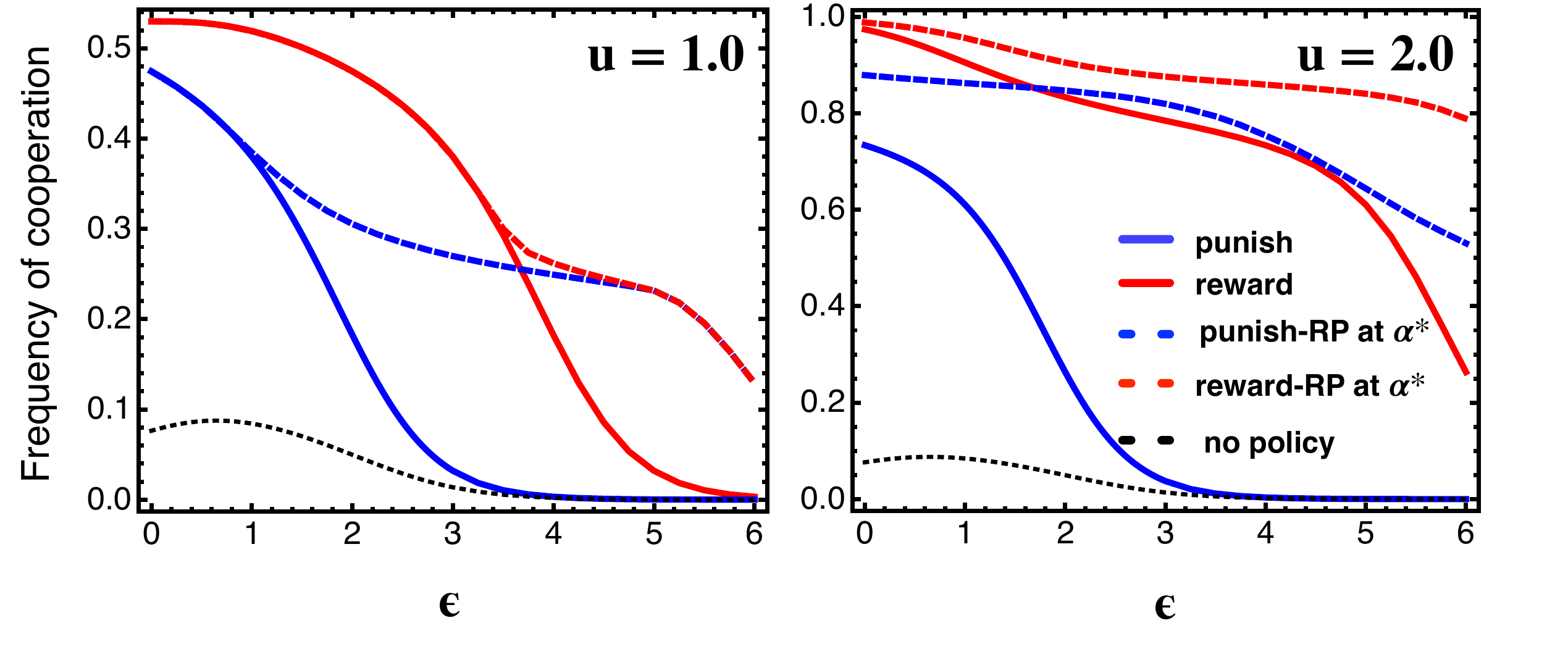}
    \caption{\textbf{Rewarding participation promotes higher levels of  cooperation given a sufficient budget ($u$)}. We compare pure reward and punishment (solid red and blue lines) and when an optimal fraction ($\alpha^\star$) of the budget is used for rewarding participation (red and blue dashed lines). When $u$ is small (left panel), the improvement occurs  when $\epsilon$ is sufficiently large. When $u$ is large (right panel), participation can  provide a large improvement for the whole range of  $\epsilon$.  The improvement obtained through rewarding participation is more significant for punishment than for reward, for both small and large $u$.   Other parameters: $N = 100$,  $R = 1, \ S = -1, \ T = 2, \ P = 0$.   }
    \label{fig:cooperation-at-optimal-alpha-vs-pure-incentives-vary-eps}
\end{figure}

\subsubsection{Incentives with rewarding participation}
We show in Figure \ref{fig:reward-vs-punishment-with-reward-participation-cooperation}  the frequency of cooperation obtained through reward or punishment when a fraction $\alpha$ of the per capita budget $u$ is used for rewarding those who agree to participate. %\ty{, regardless of whether an agreement is formed eventually.} {[We show in SI that if rewarding participation is only applied when an agreement is formed, it's detrimental to cooperation (?) because ...]}.   
Note that the scenario where $\alpha = 0$ reproduces the above described pure punishment and pure reward scenarios. For  both types of incentive, 
reward of participation in a prior commitment  can improve  overall cooperation, especially with a larger available budget for incentive supply (compare top and bottom rows). We observe that, the larger the cost of commitment participation $\epsilon$, the larger the fraction $\alpha$ of the budget should be used for rewarding participation in order for that population to reach the highest overall frequency of   cooperation.  Only when both $u$ and $\epsilon$ are sufficiently small,   {it is} better off not rewarding participation (see $\epsilon = 0.5$ and $u = 1.0$).

With $\alpha^\star$ representing the value of $\alpha$ leading to the highest level of cooperation,  
  Figure \ref{fig:cooperation-at-optimal-alpha-vs-pure-incentives-vary-eps}  compares the level of cooperation  obtained at $\alpha =  \alpha^\star$,  that when participation is not incentivised (i.e. $\alpha = 0$), and that when no policy is in place. We observe that when $u$ is small ($u = 1.0$), reward of participation leads to an improvement  when $\epsilon$ is sufficiently large. When $u$ is larger ($u = 2.0$), a significant  improvement is observed  for the whole range of  $\epsilon$ being considered. In addition, we find that the improvement obtained through the reward of participation is greater  for punishment than for reward, for both small and large $u$. This observation is in line with  the ESS analysis above.  

%Players need a push when participation is costly. 

These notable results can be explained by looking at the frequency of strategies as a function of $\alpha$, in  Figure S4 in \textbf{SI}.  As $\alpha$ approaches $\alpha^\star$, for both reward and punishment, the frequency of NDD decreases and those of ACD and ACC increase. This increase is more significant for punishment than for reward. When $\alpha >  \alpha^\star$, ADD frequency starts to increase quickly and becomes dominant in the population since the remaining budget for incentivising commitment-compliant behaviours becomes insufficient. Figure S6 in \textbf{SI} shows that this observation is robust for other values of $\epsilon$.

% \begin{figure}
%     \centering
%     \includegraphics[width=\linewidth]{panel_vary_noise.pdf}
%     \caption{\textbf{{ADD A FIGURE FOR VARYING ALPHA FOR DIFFERENT VALUES OF OMEGA; U = 2, EPS = 0.5; ALSO A CONTOUR PLOT FOR VARYING NOISE AND ALPHA?}Commitment compliance and cooperation prevails in the presence of noise.} Depicted are the frequency of  strategies and the total level of cooperation for varying  the error probability of  decision making at the commitment  stage (${\chi}$).  We consider both pure punishment (left column) and pure reward (right column), for different values of $u$ (top row, $u = 1$; bottom row, $u = 2$).  The frequency of commitment-compliant strategy (ACD)  benefits significantly from having some noise since it is now risk-dominant against ACC, which is not the case in absence of noise.  Other parameters: population size $N = 100$, $\alpha = 0$, $\beta = 0.1$, $\epsilon= 0.5$,  $R = 1, \ S = -1, \ T = 2, \ P = 0$.  }
%     \label{fig:panel_noise}
% \end{figure}

\begin{figure}
    \centering
    \includegraphics[width=\linewidth]{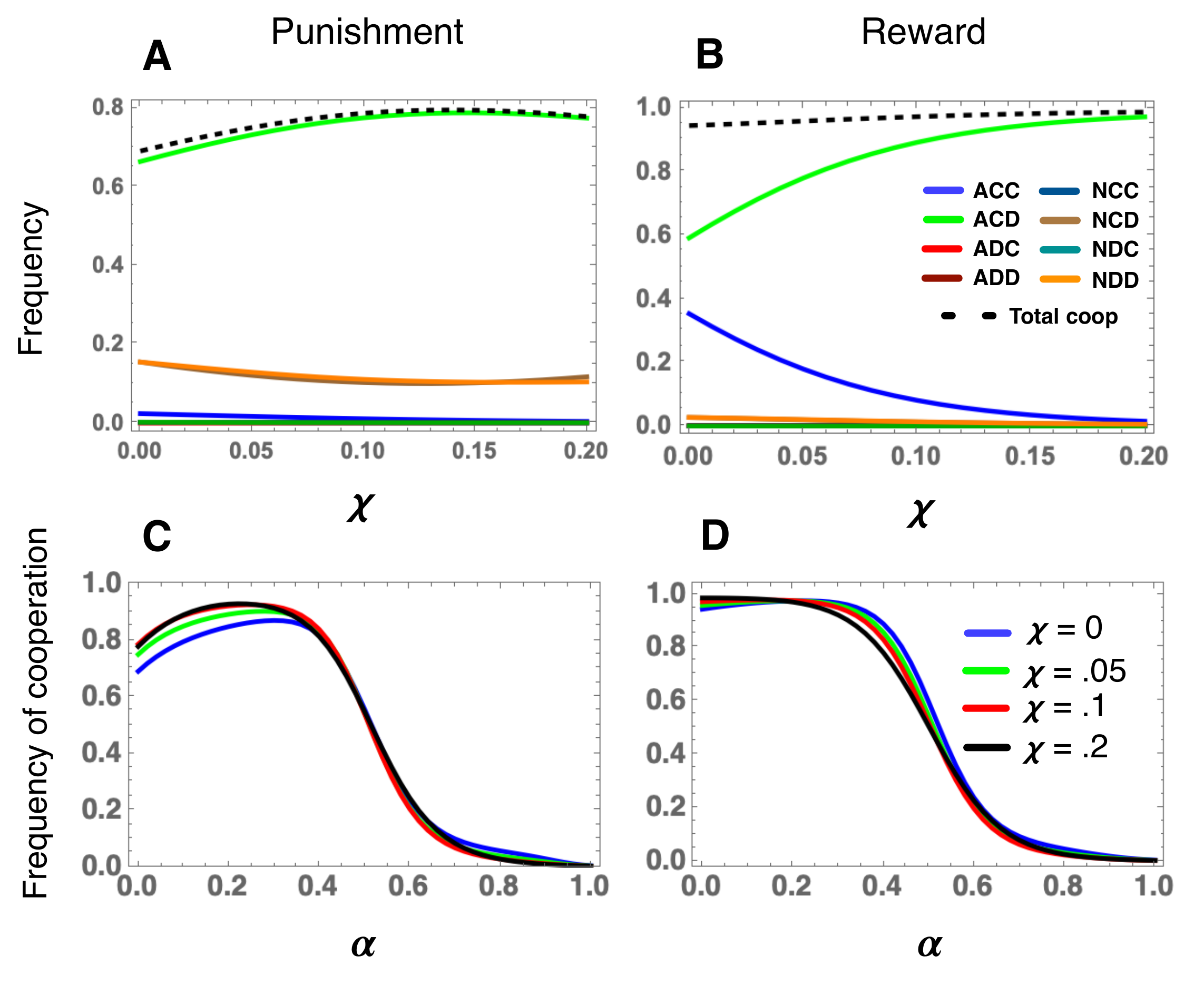}
    \caption{\textbf{Commitment compliance and cooperation prevails in the presence of noise.} Panels A, B: Frequency of  strategies and the total level of cooperation for varying  the noise probability at   the commitment  stage (${\chi}$), for pure punishment (A) and pure reward (B).  The frequency of commitment-compliant strategy (ACD)  increases in the presence of noise since it is now risk-dominant against ACC, which is not the case in absence of noise.  Panels C, D: Frequency of cooperation as a function of the fraction of the budget for rewarding participation ($\alpha$), for different noise probabilities, for punishment (C) and reward (D). A lower optimal value of $\alpha$ is observed for a greater ${\chi}$.  Other parameters: population size $N = 100$, $\alpha = 0$, $\epsilon= 0.5$, $u = 2$, $R = 1, \ S = -1, \ T = 2, \ P = 0$.  }
    \label{fig:panel_noise}
\end{figure}

\subsubsection{Non-negligible error  in  commitment participation} 
%In the ESS analysis, we have seen that when noise  in the decision making at the commitment stage is not negligible, the outcome is significantly changed since it breaks the tie between ACD and ACC. 

We study how  non-negligible error  when deciding whether to participate in a commitment (i.e., an AXY player would refuse to commit and act in the same way as a NXY player, and vice versa)  impacts the evolutionary dynamics. 
When it is negligible (${\chi} = 0$),  ACD cannot be an ESS nor risk-dominant against all other strategies in the population because  it is always neutral to ACC.  Interestingly, whenever this error probability is positive,  ACD becomes always risk-dominant against ACC (see \textbf{SI}), thus making it likely to be more prevalent.
Indeed, as shown in Figure \ref{fig:panel_noise} (panels A, B)  depicting  the frequency of the strategies and the overall level of cooperation as a function of ${\chi}$, ACD becomes more frequent and ACC less frequent as ${\chi}$ increases. The impact is more significant in case of reward than punishment. It  leads to a slight increase in the overall level of cooperation in both cases. These observations are robust for other values of $\epsilon$ and $u$ (see \textbf{SI}, Figure S7).

 Figure \ref{fig:panel_noise} shows that rewarding participation in a commitment can improve  overall cooperation in the presence of error, especially in the case of punishment. For a larger error probability, we observe a smaller value of  $\alpha^\star$ that leads to the highest level of cooperation, for both types of incentive. Moreover, we  observe a slight increase in cooperation at $\alpha^\star$  for increasing ${\chi}$. These observations are robust for other values of $\epsilon$  (see \textbf{SI}, Figure S8).
% \subsection{Extensions, robustness, other analyses }

% Different values of $\beta$, benefit to cost ratio $b/c$, importance of the games (like in Sci Rep 2013 and Trust paper)

% We study the case where reward for participation is only made when the agreement is formed in SI; i.e. only when the co-player of those who are willing to participate are also willing to do so (for the commitment to be formed).  

%Noise during the interactions
\section{Discussion}
It has been suggested that human specialised capacity for commitment might have been shaped by natural selection \citep{nesse2001evolution,frank88,akdeniz2021evolution}. Arranging a commitment from all parties involved prior to an interaction can increase the chance of reaching mutual cooperation  \citep{cherry2013enforcing,chen1994effects,dannenberg2016non,sasaki2015commitment}, enabling individuals to clarify preferences or intentions from their partners before committing to a potentially costly course of actions \citep{chen1994effects,han2015synergy,tomasello2005understanding,sterelny2012evolved}. Since individuals can decide whether not to honour an adopted commitment---there being abundant evidence  of commitment breaching in both controlled experiments and real-world scenarios \citep{nesse2001evolution,dannenberg2016non,kerr1997still,nguyen2019}---it is important to understand what mechanisms, such as positive and negative incentives, are more efficient at ensuring the cooperation-promoting benefit provided by commitments. 

%Theoretical models of commitments have shown that this is indeed the case in    social dilemma settings  \citep{han2013good,han2015synergy,Han:2014tl,bianca2021coordination_agreement,barrett2016coordination,sasaki2015commitment}.   Corresponding models, however, is limited in identifying which type of i, wherein compensation can always be enforced from those who dishonour an adopted commitment. 

%While the assumption of  regimented commitments is an useful idealisation, it can be too strict in many applications.  Commitment violators might refuse to compensate, are not capable of compensating/paying fine, or  even attempt to escape enforcement. Relaxing this assumption,

Herein, we have comparatively explored  institutional punishment of  commitment violators and reward of commitment fulfillers as potential  mechanisms to enhance commitment compliance and thus the overall cooperation in the population. We have shown that, given the same, sufficiently high, per capita budget for supplying incentives, reward results in  a higher level of commitment compliance and cooperation than punishment.  
{Reward can ensure commitment compliance to be evolutionarily stable for a larger range of the commitment cost ($\epsilon$) and benefit to cost ratio ($b/c$).}
{This observation has useful implications for the design of institutional mechanisms for promoting pro-social behaviour, especially when  communication is allowed to establish prior  commitments/agreements  \citep{chen1994effects,cherry2013enforcing,nesse2001evolution}. It is important to note that the observation was obtained under the assumption  that reward and punishment are equivalently cost-efficient (see   Model), i.e. given the same per-capita incentive budget $u$, they make the same impact in terms of payoff increase or reduction, respectively, to the incentive recipient. } 
{Also, it might be more costly for  institutions to provide rewards when compliant behaviour is frequent. Thus, an interesting direction is to consider how to combine reward and punishment in a cost efficient way, taking into account that they might have different levels of cost-efficiency,} as have been done in the context of social dilemmas (without considering commitment-based behaviours) \citep{chen2015first,DuongHanPROCsA2021,sasaki2012take,gois2019reward,sun2021combination,CIMPEANU2021107545}. Our results suggest that further behavioural experiments are needed to examine what incentive mechanisms are  efficient and preferred  (by people) for ensuring the commitment compliance and cooperation. Note that, similar to the theoretical modelling literature of incentives (see Introduction), there have been a body of behavioural experimental works comparing the effectiveness of reward and punishment for promoting cooperation in social dilemmas, see e.g. \citep{rand2009positive,balliet2011reward,van2014reward}. However, they have not studied reward and punishment in the presence of prior commitments.

Participating in a commitment can be quite costly and that might discourage players from joining the commitment in the first place. We hypothesised  that by spending part of the per capita budget to incentivise  participation before  interaction,   higher  levels of compliance and cooperation might be achievable. Indeed, we have shown that the larger the cost of participation ($\epsilon$) is, the greater fraction of this budget  should to be used for encouraging participation to achieve an optimal level of cooperation. This observation confirms the importance of studying incentives in models considering an explicit process of commitment formation, which has been omitted in extant models of institutional incentives \citep{sigmundinstitutions,chen2015first,wang2019exploring,DuongHanPROCsA2021,gois2019reward,garcia2019evolution}. Indeed, these works have not considered  commitment-based interactions nor incentives for encouraging participation in the interaction, nor its impact on the overall cooperation. %, even when participation  is optional \citep{sasaki2012take,key:Hauert2007}. 

In a cooperative interaction,  environmental noise is usually expected to lead to a  detrimental impact on the emergence and stability of cooperation \citep{nowak:2006bo,key:Sigmund_selfishnes,kahneman2021noise}  and thus requires additional supporting mechanisms such as apology and forgiveness \citep{martinez2015apology,mccullough2008beyond}. Surprisingly, we have shown here that the presence of some noise that causes errors when deciding to participate in a prior commitment, can stabilise  commitment compliance and cooperation. {In fact,  ACD can be an ESS and risk-dominant against all other strategies only when the error probability is non-zero (${\chi} > 0$)}. A non-negligible level of noise enables ACD to break ties with the commitment-accepting unconditional cooperators (ACC), who cooperate even when a commitment is not formed and thus can be easily exploited by non-committing defectors.  In SI, we also considered execution noise that happens during the PD game. We showed that it has insignificant effects on the evolutionary  dynamics and ESS analysis (e.g. none of the strategies can be an ESS for this type of noise).

A drawback of pro-social incentives for promoting cooperation is the possibility of antisocial reward and punishment  where defectors might punish cooperators or reward other defectors, hence hindering the evolution of cooperation \citep{Herrmann2008,Rand2011,dos2017antisocial,szolnoki2015antisocial}. 
We argue that this issue disappears when a prior  commitment is arranged since it  becomes clear what behaviour is expected from all parties involved during the interaction.  Only those who  commit to cooperate can be punished for  defection or rewarded for cooperation. 
It is not deemed justifiable to  punish defectors  or reward cooperators if they did commit in the first place. 
That is, commitments enable the freedom of choice from players, which can be important in cases where it might be contestable whether a behaviour is good \citep{HAN2022101843}, or when players might not be capable of cooperation, for example due to other commitments or lack of resources to carry it out.

{When employing institutional incentives for sustaining  pro-social behaviours, an important issue is how to set up and maintain a sufficient incentive budget. The problem of who contributes to this budget is a social dilemma in itself, and how to overcome it is a challenging  problem. Several solutions have been proposed recently,  including pool incentives with second order punishments \citep{sigmundinstitutions}, democratic decisions \citep{hilbe2014democratic} and positive and negative incentives combination \citep{gois2019reward,chen2015first}. In this work, we assume that an institution exists to facilitate the incentive providing process for promoting compliance and participation in an agreement, where those who agree to join a prior commitment contribute a participation fee ($\epsilon$) to help sustain the institution. 
Future works might consider other non-institutional mechanisms that might underline   commitment compliance and participation, such as  reputation-based mechanisms (for instance, those who dishonour an adopted commitment are assigned low reputation scores)  \citep{nowak2005evolution,perret2021evolution,okada2020review} and emotional incentives (e.g. players might feel more guilty about breaching a commitment than a mere wrongdoing without a prior commitment) \citep{Luis2017AAMAS,vanberg2008people,nesse2001evolution}, and even how these mechanisms might interplay and be combined with the institutional approach for further improvement \citep{nesse2001evolution}. 
}

Evolutionary modelling and analysis of  voluntary participation has been considered in several studies \citep{de2010freedom,mathew2009does,sasaki2012take,key:Hauert2007,sigmundinstitutions,salahshour2021evolution}, showing that cooperation can evolve even in one-shot cooperation dilemmas if players have the option to  opt out.
However, these works did not  consider strategies  conditioned  on the formation of a commitment, nor incentives for encouraging the participation in it. Typically only a subset of unconditional strategies  were considered, including  cooperators,  defectors and non-participants. Here we have examined a full set of strategies (Table \ref{table:eight-strategies}). In contrast to these studies, we have shown that the evolutionary stable strategies often exhibit  behaviours conditional on the formation of commitment being formed, e.g., ACD and NDC are ESS (see again Figure \ref{fig:ESS_analysis_overall}). Given this crucial limitation of previous works, our model here provides a more complete picture of how prior commitments such as formal and informal contracts and agreements, provides an efficient mechanism for promoting the evolution of cooperation. 

In short, we have analysed here different forms of institutional incentive for promoting participation and compliance in interactions with a prior commitment formation. Our results have shown that, in this setting, using incentives for ensuring participation is as important as for enhancing compliance.   
 
% We consider rewarding participation in the paper. Sanctioning non-participation is an option but it can be considered unethical/unjustifiable given freedom of choice rights (any refs). 

%Commitments have been widely studied in the context  of  multi-agent and autonomous computerised agent systems, in order to ensure high levels of cooperation among agents Wooldridge and Jennings 1999; Winikoff 2007; MORE They have also been utilized for ensuring good behaviors in various computerised applications such as electric vehicle charging (Stein et al. 2012) and peer-to-peer sharing networks (Rzadca et al. 2015).

%\footnote{In multi party interactions, it would be different as one can envisage an agreement can be formed if a threshold number of players accepted to participate (see jaamas paper). (this is however beyond the scope of this paper). They would also be different if we consider noise in the pre-commitment stage, when players decide whether or not to join an agreement (e.g. due to traffic/internet failure/miscommunication, etc), or during the games (in or out of agreement). }.

%Future works: multi-player games, climate change, AI risks; interplay between incentives for mere cooperators and commitment compliant players 

%\newpage
\section*{Author Contributions}
T.A.H. designed the research, developed the model, implemented software,  carried out the analysis and wrote the manuscript.
\section*{Acknowledgements}
The author is grateful for the useful discussion and comments from Luis Moniz Pereira, Jeff White and Paolo Bova on early versions of this manuscript.   

\section*{Data Accessibility}
This work does not contain any data.

\section*{Funding Statement.}
This research is supported by a Leverhulme Research Fellowship   ``Incentives for Commitment Compliance" (RF-2020-603/9).

%\bibliographystyle{apalike}
%\bibliography{bibliography} 

\begin{thebibliography}{}

\bibitem[Akdeniz and van Veelen, 2021]{akdeniz2021evolution}
Akdeniz, A. and van Veelen, M. (2021).
\newblock The evolution of morality and the role of commitment.
\newblock {\em Evolutionary Human Sciences}, pages 1--53.

\bibitem[Balliet, 2010]{balliet2010communication}
Balliet, D. (2010).
\newblock Communication and cooperation in social dilemmas: A meta-analytic
  review.
\newblock {\em Journal of Conflict Resolution}, 54(1):39--57.

\bibitem[Balliet et~al., 2011]{balliet2011reward}
Balliet, D., Mulder, L.~B., and Van~Lange, P.~A. (2011).
\newblock Reward, punishment, and cooperation: a meta-analysis.
\newblock {\em Psychological bulletin}, 137(4):594.

\bibitem[Barrett and Stavins, 2003]{barrett2003increasing}
Barrett, S. and Stavins, R. (2003).
\newblock Increasing participation and compliance in international climate
  change agreements.
\newblock {\em International Environmental Agreements}, 3(4):349--376.

\bibitem[Bruni et~al., 2009]{bruni2009economic}
Bruni, M.~L., Nobilio, L., and Ugolini, C. (2009).
\newblock Economic incentives in general practice: the impact of
  pay-for-participation and pay-for-compliance programs on diabetes care.
\newblock {\em Health policy}, 90(2-3):140--148.

\bibitem[Chen et~al., 2015]{chen2015first}
Chen, X., Sasaki, T., Br{\"a}nnstr{\"o}m, {\AA}., and Dieckmann, U. (2015).
\newblock First carrot, then stick: how the adaptive hybridization of
  incentives promotes cooperation.
\newblock {\em Journal of The Royal Society Interface}, 12(102):20140935.

\bibitem[Chen and Komorita, 1994]{chen1994effects}
Chen, X.-P. and Komorita, S.~S. (1994).
\newblock The effects of communication and commitment in a public goods social
  dilemma.
\newblock {\em Organizational Behavior and Human Decision Processes},
  60(3):367--386.

\bibitem[Cherry and McEvoy, 2013]{cherry2013enforcing}
Cherry, T.~L. and McEvoy, D.~M. (2013).
\newblock Enforcing compliance with environmental agreements in the absence of
  strong institutions: An experimental analysis.
\newblock {\em Environmental and Resource Economics}, 54(1):63--77.

\bibitem[Cimpeanu et~al., 2021]{CIMPEANU2021107545}
Cimpeanu, T., Perret, C., and Han, T.~A. (2021).
\newblock Cost-efficient interventions for promoting fairness in the ultimatum
  game.
\newblock {\em Knowledge-Based Systems}, 233:107545.

\bibitem[Coombs, 1973]{coombs1973reparameterization}
Coombs, C.~H. (1973).
\newblock A reparameterization of the prisoner's dilemma game.
\newblock {\em Behavioral Science}, 18(6):424--428.

\bibitem[Dannenberg, 2016]{dannenberg2016non}
Dannenberg, A. (2016).
\newblock Non-binding agreements in public goods experiments.
\newblock {\em Oxford Economic Papers}, 68(1):279--300.

\bibitem[De~Silva et~al., 2010]{de2010freedom}
De~Silva, H., Hauert, C., Traulsen, A., and Sigmund, K. (2010).
\newblock Freedom, enforcement, and the social dilemma of strong altruism.
\newblock {\em Journal of Evolutionary Economics}, 20(2):203--217.

\bibitem[Dos~Santos and Pe{\~n}a, 2017]{dos2017antisocial}
Dos~Santos, M. and Pe{\~n}a, J. (2017).
\newblock Antisocial rewarding in structured populations.
\newblock {\em Scientific Reports}, 7(1):1--14.

\bibitem[Duong and Han, 2021]{DuongHanPROCsA2021}
Duong, M.~H. and Han, T.~A. (2021).
\newblock Cost efficiency of institutional incentives for promoting cooperation
  in finite populations.
\newblock {\em Proceedings of the Royal Society A}, 477(2254):20210568.

\bibitem[Frank, 1988]{frank88}
Frank, R.~H. (1988).
\newblock {\em Passions {W}ithin {R}eason: {T}he {S}trategic {R}ole of the
  {E}motions}.
\newblock Norton and Company.

\bibitem[Garc{\'\i}a and Traulsen, 2019]{garcia2019evolution}
Garc{\'\i}a, J. and Traulsen, A. (2019).
\newblock Evolution of coordinated punishment to enforce cooperation from an
  unbiased strategy space.
\newblock {\em Journal of the Royal Society Interface}, 16(156):20190127.

\bibitem[G{\'o}is et~al., 2019]{gois2019reward}
G{\'o}is, A.~R., Santos, F.~P., Pacheco, J.~M., and Santos, F.~C. (2019).
\newblock Reward and punishment in climate change dilemmas.
\newblock {\em Sci. Rep.}, 9(1):1--9.

\bibitem[Han, 2013]{HanBook2013}
Han, T.~A. (2013).
\newblock {\em Intention Recognition, Commitments and Their Roles in the
  Evolution of Cooperation: From Artificial Intelligence Techniques to
  Evolutionary Game Theory Models}, volume~9.
\newblock Springer SAPERE series.

\bibitem[Han and Lenaerts, 2016]{hanTom2016synergy}
Han, T.~A. and Lenaerts, T. (2016).
\newblock A synergy of costly punishment and commitment in cooperation
  dilemmas.
\newblock {\em Adaptive Behavior}, 24(4):237--248.

\bibitem[Han et~al., 2022]{HAN2022101843}
Han, T.~A., Lenaerts, T., Santos, F.~C., and Pereira, L.~M. (2022).
\newblock Voluntary safety commitments provide an escape from over-regulation
  in ai development.
\newblock {\em Technology in Society}, 68:101843.

\bibitem[Han et~al., 2015a]{Han:2014tl}
Han, T.~A., Pereira, L.~M., and Lenaerts, T. (2015a).
\newblock {Avoiding or Restricting Defectors in Public Goods Games?}
\newblock {\em J. Royal Soc Interface}, 12(103):20141203.

\bibitem[Han et~al., 2017]{HanJaamas2016}
Han, T.~A., Pereira, L.~M., and Lenaerts, T. (2017).
\newblock Evolution of commitment and level of participation in public goods
  games.
\newblock {\em Autonomous Agents and Multi-Agent Systems}, pages 1--23.

\bibitem[Han et~al., 2013]{han2013good}
Han, T.~A., Pereira, L.~M., Santos, F.~C., and Lenaerts, T. (2013).
\newblock Good agreements make good friends.
\newblock {\em Scientific reports}, 3(2695).

\bibitem[Han et~al., 2015b]{han2015synergy}
Han, T.~A., Santos, F.~C., Lenaerts, T., and Pereira, L.~M. (2015b).
\newblock Synergy between intention recognition and commitments in cooperation
  dilemmas.
\newblock {\em Scientific reports}, 5(9312).

\bibitem[Hauert et~al., 2007]{key:Hauert2007}
Hauert, C., Traulsen, A., Brandt, H., Nowak, M.~A., and Sigmund, K. (2007).
\newblock Via freedom to coercion: The emergence of costly punishment.
\newblock {\em Science}, 316:1905--1907.

\bibitem[Heidar, 2006]{heidar2006party}
Heidar, K. (2006).
\newblock Party membership and participation.
\newblock {\em Handbook of party politics}, pages 301--315.

\bibitem[Herrmann et~al., 2008]{Herrmann2008}
Herrmann, B., Th\"{o}ni, C., and G\"{a}chter, S. (2008).
\newblock {Antisocial Punishment Across Societies}.
\newblock {\em Science}, 319(5868):1362--1367.

\bibitem[Hilbe et~al., 2014]{hilbe2014democratic}
Hilbe, C., Traulsen, A., R{\"o}hl, T., and Milinski, M. (2014).
\newblock Democratic decisions establish stable authorities that overcome the
  paradox of second-order punishment.
\newblock {\em PNAS}, 111(2):752--756.

\bibitem[Hofbauer and Sigmund, 1998]{key:Hofbauer1998}
Hofbauer, J. and Sigmund, K. (1998).
\newblock {\em Evolutionary Games and Population Dynamics}.
\newblock Cambridge University Press.

\bibitem[Imhof et~al., 2005]{key:imhof2005}
Imhof, L.~A., Fudenberg, D., and Nowak, M.~A. (2005).
\newblock Evolutionary cycles of cooperation and defection.
\newblock {\em Proc. Natl. Acad. Sci. U.S.A.}, 102:10797--10800.

\bibitem[Irons, 2001]{IronsChapter2001}
Irons, W. (2001).
\newblock Religion as a hard-to-fake sign of commitment.
\newblock In Nesse, R.~M., editor, {\em Evolution and the capacity for
  commitment}, pages 292--309. New York: Russell Sage.

\bibitem[Johnson and Bering, 2006]{johnson2006hand}
Johnson, D. and Bering, J. (2006).
\newblock Hand of god, mind of man: Punishment and cognition in the evolution
  of cooperation.
\newblock {\em Evolutionary psychology}, 4(1):147470490600400119.

\bibitem[Kahneman et~al., 2021]{kahneman2021noise}
Kahneman, D., Sibony, O., and Sunstein, C.~R. (2021).
\newblock {\em Noise: a flaw in human judgment}.
\newblock Little, Brown.

\bibitem[Kerr et~al., 1997]{kerr1997still}
Kerr, N.~L., Garst, J., Lewandowski, D.~A., and Harris, S.~E. (1997).
\newblock That still, small voice: Commitment to cooperate as an internalized
  versus a social norm.
\newblock {\em Personality and social psychology Bulletin}, 23(12):1300--1311.

\bibitem[Martinez-Vaquero et~al., 2015]{martinez2015apology}
Martinez-Vaquero, L.~A., Han, T.~A., Pereira, L.~M., and Lenaerts, T. (2015).
\newblock Apology and forgiveness evolve to resolve failures in cooperative
  agreements.
\newblock {\em Scientific reports}, 5(10639).

\bibitem[Mathew and Boyd, 2009]{mathew2009does}
Mathew, S. and Boyd, R. (2009).
\newblock When does optional participation allow the evolution of cooperation?
\newblock {\em Proceedings of the Royal Society B: Biological Sciences},
  276(1659):1167--1174.

\bibitem[Maynard-Smith, 1982]{maynard-smith:1982to}
Maynard-Smith, J. (1982).
\newblock {\em Evolution and the Theory of Games}.
\newblock Cambridge University Press, Cambridge.

\bibitem[McCullough, 2008]{mccullough2008beyond}
McCullough, M. (2008).
\newblock {\em Beyond revenge: The evolution of the forgiveness instinct}.
\newblock John Wiley \& Sons.

\bibitem[Nesse, 2001]{nesse2001evolution}
Nesse, R.~M. (2001).
\newblock {\em Evolution and the capacity for commitment}.
\newblock Foundation series on trust. Russell Sage.

\bibitem[Nguyen et~al., 2019]{nguyen2019}
Nguyen, H.~K., Chiong, R., Chica, M., Middleton, R., and Thi Kim~Pham, D.
  (2019).
\newblock Contract farming in the mekong delta's rice supply chain: Insights
  from an agent-based modeling study.
\newblock {\em Journal of Artificial Societies and Social Simulation}, 22(3):1.

\bibitem[Nowak, 2006]{nowak:2006bo}
Nowak, M.~A. (2006).
\newblock {\em Evolutionary Dynamics: Exploring the Equations of Life}.
\newblock Harvard University Press, Cambridge, MA.

\bibitem[Nowak et~al., 2004]{key:novaknature2004}
Nowak, M.~A., Sasaki, A., Taylor, C., and Fudenberg, D. (2004).
\newblock Emergence of cooperation and evolutionary stability in finite
  populations.
\newblock {\em Nature}, 428:646--650.

\bibitem[Nowak and Sigmund, 2005a]{nowak:2005:nature}
Nowak, M.~A. and Sigmund, K. (2005a).
\newblock Evolution of indirect reciprocity.
\newblock {\em Nature}, 437(1291-1298).

\bibitem[Nowak and Sigmund, 2005b]{nowak2005evolution}
Nowak, M.~A. and Sigmund, K. (2005b).
\newblock Evolution of indirect reciprocity.
\newblock {\em Nature}, 437(7063):1291--1298.

\bibitem[Okada, 2020]{okada2020review}
Okada, I. (2020).
\newblock A review of theoretical studies on indirect reciprocity.
\newblock {\em Games}, 11(3):27.

\bibitem[Ostrom, 2005]{ostrom2009understanding}
Ostrom, E. (2005).
\newblock {\em Understanding institutional diversity}.
\newblock Princeton university press.

\bibitem[Otto and Day, 2007]{Otto2007AEvolution}
Otto, S.~P. and Day, T. (2007).
\newblock {\em {A biologist's guide to mathematical modeling in ecology and
  evolution}}, volume~6.
\newblock Princeton University Press, Princeton, NJ.

\bibitem[Pereira et~al., 2017]{Luis2017AAMAS}
Pereira, L.~M., Lenaerts, T., Martinez-Vaquero, L.~A., and Han, T.~A. (2017).
\newblock Social manifestation of guilt leads to stable cooperation in
  multi-agent systems.
\newblock In {\em AAMAS}, pages 1422--1430.

\bibitem[Perret et~al., 2021]{perret2021evolution}
Perret, C., Krellner, M., and Han, T.~A. (2021).
\newblock The evolution of moral rules in a model of indirect reciprocity with
  private assessment.
\newblock {\em Scientific Reports}, 11(1):1--10.

\bibitem[Rand et~al., 2009]{rand2009positive}
Rand, D.~G., Dreber, A., Ellingsen, T., Fudenberg, D., and Nowak, M.~A. (2009).
\newblock Positive interactions promote public cooperation.
\newblock {\em Science}, 325(5945):1272--1275.

\bibitem[Rand and Nowak, 2011]{Rand2011}
Rand, D.~G. and Nowak, M.~A. (2011).
\newblock The evolution of antisocial punishment in optional public goods
  games.
\newblock {\em Nature Communications}, 2:434.

\bibitem[Rand et~al., 2013]{randUltimatum}
Rand, D.~G., Tarnita, C.~E., Ohtsuki, H., and Nowak, M.~A. (2013).
\newblock Evolution of fairness in the one-shot anonymous ultimatum game.
\newblock {\em Proc. Natl. Acad. Sci. USA}, 110:2581--2586.

\bibitem[Salahshour, 2021]{salahshour2021evolution}
Salahshour, M. (2021).
\newblock Evolution of cooperation in costly institutions exhibits red queen
  and black queen dynamics in heterogeneous public goods.
\newblock {\em Communications biology}, 4(1):1--10.

\bibitem[Sasaki et~al., 2012]{sasaki2012take}
Sasaki, T., Br{\"a}nnstr{\"o}m, {\AA}., Dieckmann, U., and Sigmund, K. (2012).
\newblock The take-it-or-leave-it option allows small penalties to overcome
  social dilemmas.
\newblock {\em Proceedings of the National Academy of Sciences},
  109(4):1165--1169.

\bibitem[Sasaki et~al., 2015]{sasaki2015commitment}
Sasaki, T., Okada, I., Uchida, S., and Chen, X. (2015).
\newblock Commitment to cooperation and peer punishment: Its evolution.
\newblock {\em Games}, 6(4):574--587.

\bibitem[Sehgal, 2021]{sehgal2021impact}
Sehgal, N.~K. (2021).
\newblock Impact of vax-a-million lottery on covid-19 vaccination rates in
  ohio.
\newblock {\em The American Journal of Medicine}.

\bibitem[Sigmund, 2010]{key:Sigmund_selfishnes}
Sigmund, K. (2010).
\newblock {\em The Calculus of Selfishness}.
\newblock Princeton University Press.

\bibitem[Sigmund et~al., 2010]{sigmundinstitutions}
Sigmund, K., Silva, H.~D., Traulsen, A., and Hauert, C. (2010).
\newblock Social learning promotes institutions for governing the commons.
\newblock {\em Nature}, 466:7308.

\bibitem[Singh, 2013]{singh2013norms}
Singh, M.~P. (2013).
\newblock Norms as a basis for governing sociotechnical systems.
\newblock {\em ACM Transactions on Intelligent Systems and Technology (TIST)},
  5(1):21.

\bibitem[Sterelny, 2012]{sterelny2012evolved}
Sterelny, K. (2012).
\newblock {\em The evolved apprentice}.
\newblock MIT Press.

\bibitem[Sun et~al., 2021]{sun2021combination}
Sun, W., Liu, L., Chen, X., Szolnoki, A., and Vasconcelos, V.~V. (2021).
\newblock Combination of institutional incentives for cooperative governance of
  risky commons.
\newblock {\em Iscience}, 24(8):102844.

\bibitem[Szab{\'o} and F{\'a}th, 2007]{Szabo2007}
Szab{\'o}, G. and F{\'a}th, G. (2007).
\newblock Evolutionary games on graphs.
\newblock {\em Phys Rep}, 97-216(4-6).

\bibitem[Szolnoki and Perc, 2015]{szolnoki2015antisocial}
Szolnoki, A. and Perc, M. (2015).
\newblock Antisocial pool rewarding does not deter public cooperation.
\newblock {\em Proceedings of the Royal Society B: Biological Sciences},
  282(1816):20151975.

\bibitem[Tappin et~al., 2015]{tappin2015financial}
Tappin, D., Bauld, L., Purves, D., Boyd, K., Sinclair, L., MacAskill, S.,
  McKell, J., Friel, B., McConnachie, A., De~Caestecker, L., et~al. (2015).
\newblock Financial incentives for smoking cessation in pregnancy: randomised
  controlled trial.
\newblock {\em Bmj}, 350.

\bibitem[Tomasello et~al., 2005]{tomasello2005understanding}
Tomasello, M., Carpenter, M., Call, J., Behne, T., and Moll, H. (2005).
\newblock Understanding and sharing intentions: The origins of cultural
  cognition.
\newblock {\em Behavioral and brain sciences}, 28(05):675--691.

\bibitem[Traulsen et~al., 2006]{traulsen2006}
Traulsen, A., Nowak, M.~A., and Pacheco, J.~M. (2006).
\newblock Stochastic dynamics of invasion and fixation.
\newblock {\em Phys. Rev. E}, 74:11909.

\bibitem[Van~Lange et~al., 2014]{van2014reward}
Van~Lange, P.~A., Rockenbach, B., and Yamagishi, T. (2014).
\newblock {\em Reward and punishment in social dilemmas}.
\newblock Oxford University Press.

\bibitem[Vanberg, 2008]{vanberg2008people}
Vanberg, C. (2008).
\newblock Why do people keep their promises? an experimental test of two
  explanations 1.
\newblock {\em Econometrica}, 76(6):1467--1480.

\bibitem[Wang et~al., 2019]{wang2019exploring}
Wang, S., Chen, X., and Szolnoki, A. (2019).
\newblock Exploring optimal institutional incentives for public cooperation.
\newblock {\em Communications in Nonlinear Science and Numerical Simulation},
  79:104914.

\bibitem[Zisis et~al., 2015]{zisisSciRep2015}
Zisis, I., Guida, S.~D., Han, T.~A., Kirchsteiger, G., and Lenaerts, T. (2015).
\newblock Generosity motivated by acceptance - evolutionary analysis of an
  anticipation games.
\newblock {\em Scientific reports}, 5(18076).

\bibitem[Zumbansen, 2007]{zumbansen2007law}
Zumbansen, P. (2007).
\newblock The law of society: governance through contract.
\newblock {\em Indiana Journal of Global Legal Studies}, 14(2):191--233.

\end{thebibliography}

\end{document}